\documentclass[oneside]{amsart}

\usepackage{geometry}
\usepackage{amssymb}
\usepackage[mathscr]{euscript}
\usepackage{amsmath}
\usepackage{mathtools}
\usepackage{enumitem} 
\usepackage{mathrsfs}
\usepackage{lipsum}
\usepackage{amsfonts}
\usepackage{graphicx}
\usepackage{epstopdf}
\usepackage{algorithmic}
\usepackage{hyperref}

\usepackage{tikz}
\usetikzlibrary{arrows,calc}

\tikzstyle{arc}=[->, shorten <=1pt, shorten >=1pt, >=stealth, line width=1pt]
\tikzstyle{edge}=[shorten <=1pt, shorten >=1pt, >=stealth, line width=1pt]
\tikzstyle{vertex}=[circle, fill=black, draw, minimum size=4pt, inner sep=0pt, outer sep=0pt]

\newcommand{\cD}{\mathcal D}
\newcommand{\cK}{\mathcal{K}}
\newcommand{\cF}{\mathcal{F}}
\newcommand{\HF}{H_\mathcal{F}}
\newcommand{\DF}{D_\mathcal{F}}
\newcommand{\Pol}{\operatorname{Pol}}
\newcommand{\CSP}{\operatorname{CSP}}
\newcommand{\mino}{\operatorname{minority}}
\newcommand{\majo}{\operatorname{majority}}
\newcommand{\sP}{\mathscr{P}}
\newcommand{\sC}{\mathscr{C}}
\newcommand{\sD}{\mathscr{D}}
\newcommand{\sG}{\mathscr{G}}
\newcommand{\CDF}{\sC^{\DF}}
\newcommand{\CDFinj}{\sC^{\DF,\mathrm{inj}}}

\newcommand{\CAinj}{\sC^{\bA,\mathrm{inj}}}
\newcommand{\Aut}{\operatorname{Aut}}
\newcommand{\End}{\operatorname{End}}
\newcommand{\xto}{\leftarrow}
\newcommand{\Sym}{\operatorname{Sym}}
\newcommand{\bA}{\mathbb A}
\newcommand{\bB}{\mathbb B}
\newcommand{\bC}{\mathbb C}
\newcommand{\bF}{\mathbb F}
\newcommand{\bQ}{\mathbb Q}
\newcommand{\bH}{\mathbb H}
\newcommand{\bT}{\mathbb T}
\newcommand{\bX}{\mathbb X}
\newcommand{\age}{\operatorname{Age}}

\theoremstyle{plain}
\newtheorem{theorem}{Theorem}
\newtheorem{lemma}[theorem]{Lemma}
\newtheorem{corollary}[theorem]{Corollary}
\newtheorem{proposition}[theorem]{Proposition}

\theoremstyle{definition}
\newtheorem{definition}[theorem]{Definition}

\title[The graph orientation problem dichotomy for forbidden tournaments]{An algebraic proof of the dichotomy for graph orientation problems with forbidden tournaments}
\thanks{Funded by the European Union (ERC, POCOCOP, 101071674). Views and opinions expressed are however those of the author(s) only and do not necessarily reflect those of the European Union or the European Research Council Executive Agency. Neither the European Union nor the granting authority can be held responsible for them. This research was funded in whole or in part by the Austrian Science Fund (FWF) [I 5948]. For the purpose of Open Access, the authors have applied a CC BY public copyright licence to any Author Accepted Manuscript (AAM) version arising from this submission.}

\author{Roman Feller}
\author{Michael Pinsker}

\begin{document}

\begin{abstract}
For a set $\cF$ of finite tournaments, the $\cF$-free orientation problem is the problem of deciding if a given finite undirected graph can be oriented in such a way that the resulting oriented graph does not contain any member of $\cF$. Using the theory of smooth approximations, we give a new shorter proof of the complexity dichotomy for such problems obtained recently by  Bodirsky and Guzm\'{a}n-Pro. In fact, our approach  yields a complexity dichotomy for a considerably larger class of computational problems where one is given an undirected graph along with additional local constraints on the allowed orientations. Moreover, the border between tractable and hard problems is also described by a decidable algebraic condition.
\end{abstract}

\maketitle

\section{Introduction}
\subsection{Graph orientation problems} 
For a set $\cF$ of oriented graphs, the \emph{$\cF$-free orientation problem} is the following decision problem: given a finite undirected graph $G$, can  the edges of $G$ be oriented in such a way that the resulting oriented graph does not contain any member of $\cF$ as an induced subgraph? 
Two natural questions arise in this context. On the one hand, one might aim for a structural classification of all graphs admitting an $\cF$-free orientation. 
Early results in this direction are due to Skrien~\cite{Skrien}, who provides a characterisation of those graphs with an $\cF$-free orientation, where $\cF$ ranges over all possible combinations of oriented paths on three vertices. 
Another classical result, the Gallai–Hasse–Roy–Vitaver Theorem~\cite{Gallai,Hasse,Roy,Vitaver}, implies that for every $k\geq 1$ there is a set of oriented graphs $\cF_k$ such that a graph is $k$-colourable if and only if it admits an $\cF_k$-free orientation. On the other hand, one can investigate the  computational complexity of the $\cF$-free orientation problem. In the case where $\cF$ consists of oriented paths on three vertices only, the $\cF$-free orientation problem admits a polynomial-time reduction to 2-SAT~\cite{BangJensen}, and hence is tractable. In the general case, however, the $\cF$-free orientation problem might be equivalent to an NP-complete vertex coloring problem, as mentioned above. 
Recently, Bodirsky and Guzm\'an-Pro have shown that if $\cF$ consists of tournaments, then the $\cF$-free orientation problem is  either solvable in polynomial time  or NP-complete~\cite{forbiddenTournaments}, thus establishing a complexity dichotomy. 

The present paper provides a new  shorter proof of this complexity dichotomy using the novel algebraic theory of smooth approximations first announced in~\cite{MottetPinskerSmoothConf}. It thereby continues the program initiated in that work of unifying all proofs of complexity dichotomies for problems within the Bodirsky-Pinsker conjecture (see Section~\ref{sect:relatedwork} below), with the ultimate goal of removing more and more ad hoc arguments. Indeed, contrary to the proof of Bodirsky and Guzm\'{a}n-Pro that has little conceptual overlap with that of other complexity dichotomies, our proof lines up smoothly with other proofs using the theory of smooth approximations as  presented in~\cite{Hypergraphs,RELATIONALWIDTHCOLLAPSES,MottetPinskerSmooth}. It thus further contributes to a clearer picture of the general borderline between tractability and hardness, and makes the result amenable to  generalisations.

In particular, such generalisations include problems where we are not necessarily given a graph, but more generally a list of constraints on  vertices given by quantifier-free formulas over the language of graphs, in the spirit of Graph-SAT problems classified  in~\cite{BodPin-Schaefer-both}. In this direction, we achieve in the present paper a complexity classification for problems where we are given a graph enhanced with additional local constraints which prescribe some labelled cliques of the given graph to belong to a fixed set of allowed tournaments after orientation.
An example in this scope, whose computational complexity was already investigated by Bodir\-sky and Guzm\'{a}n-Pro in~\cite{forbiddenTournaments}, is the \emph{$\cF$-free orientation completion problem}. Here, one fixes a finite set of forbidden tournaments $\cF$ and is then given a finite input graph $G$ with some of its edges oriented already. The task then consists of completing this partial orientation to an $\cF$-free orientation of $G$. In this case the local constraints imposed on the input are limited to single edges, corresponding precisely to the pre-oriented edges of the input graph. 
Instead of being limited to local constraints imposed on single edges, our result allows to impose local constraints on cliques of arbitrary size. 
For example, the problem of finding an $\cF$-free orientation of an input graph $G$ such that some labelled $3$-cliques of $G$ inherit a transitive orientation while some labelled $4$-cliques inherit a cycle somewhere, lies within the scope of our result. More generally, for arbitrary $n$, we may label $n$-tuples of  vertices that induce a clique in the input graph $G$ and specify a list of allowed configurations, i.e., tournaments on $\{1,\dots,n\}$, such that after orientation any labelled tuple, say $(v_1,\dots,v_n)$ is isomorphic to an allowed configuration via the mapping $v_i \mapsto i$. This produces ``generalised graph orientation problems'' where containment in P is witnessed by structural symmetries of the \emph{template} (defined in Section~\ref{sect:CSP}) which were previously not needed in the complexity classification of graph orientation (completion) problems by Bodirsky and Guzm\'{a}n-Pro in~\cite{forbiddenTournaments}. We refer to Section~\ref{sect:examples} for a concrete example of this phenomenon. 
Potential further generalisations for which this work could well serve as a basis shall be discussed in Section~\ref{sect:outlook}.

\subsection{Constraint Satisfaction Problems (CSPs)}\label{sect:CSP} 
For a relational structure $\bA$ in some finite relational signature $\tau$, the \emph{Constraint Satisfaction Problem of $\bA$}, denoted by  $\CSP(\bA)$, is the membership problem for the class of all finite $\tau$-structures which have a homomorphism into $\bA$; the structure $\bA$ is called a \emph{template} of the CSP. 

For certain $\cF$ (including all sets of tournaments), the $\cF$-free orientation problem can be viewed as the CSP of an (in general necessarily infinite) undirected  graph with strong model-theoretic properties. We elaborate on this claim below.  This perspective provides access to the rich literature on infinite-domain constraint satisfaction.

It is easy to see that the class $\cK$ of finite $\cF$-free oriented graphs has the \emph{hereditary property}, i.e., is closed under (induced) substructures. 
If all elements of $\cF$ are connected, then $\cK$ is closed under taking disjoint unions. Consequently, there exists an infinite oriented graph $D$ whose finite substructures are precisely the elements of $\cK$ up to isomorphism. This oriented graph can be obtained simply by taking the disjoint union of all structures in $\cK$.  Let $H$ be the undirected graph reduct of $D$.\footnote{$H$ is obtained by forgetting the orientation of the edges of $D$, put differently, the edge relation of $H$ is the symmetric closure of the edge relation of $D$.} For any finite undirected graph  $G$  we then have that $G$ embeds into $H$ if and only if it is a positive instance of the $\cF$-free orientation problem. If $\cF$ consists of tournaments (or more generally, is closed under homomorphic images), this is the case if and only if $G$ has a homomorphism into $H$.  We thus have that the $\cF$-free orientation problem is the same as $\CSP(H)$. 
Finally, if $\cF$ consists of tournaments, whose underlying undirected graph is complete, then the class $\cK$ is even closed under taking unions. Equivalently, in model-theoretic jargon,  the class $\cK$ has the free amalgamation property. Consequently, the infinite oriented graph $D$ can be chosen to be \emph{homogeneous}: any finite  isomorphism between finite induced substructures of $D$ extends to an automorphism of $D$ (see Fra\"{i}ss\'{e}'s theorem in~\cite{Hodges}). 
In fact, there is only one such structure up to isomorphism which we denote by $\DF$. Accordingly, we denote the undirected graph reduct of $\DF$ by $\HF$; thus the $\cF$-free orientation problem is equal to $\CSP(\HF)$.

The structure $\HF$ has the advantage that it lies within the scope of the infinite-domain CSP dichotomy conjecture~\cite{BPP-projective-homomorphisms}, and thus the methods developed for that scope for more than a decade apply. Indeed, one of the goals of the present paper is to apply the full power of the available techniques, thus aligning the proof of the complexity classification for  $\CSP(\HF)$, and in fact  generalisations thereof, with the bigger picture both conceptually and technically. 

\subsection{The results}
We start by stating the result for $\CSP(\HF)$; however, before stating it in full detail (Theorem~\ref{sellingmainthm}), we first provide an informal description.  Our complexity dichotomy  distinguishes the following two cases. In the first case, $\HF$ can \emph{pp-construct}, i.e., encode in a certain precise way, any finite structure~\cite{wonderland}; it then follows immediately that $\CSP(\HF)$ is NP-complete. 

In the second case---excluding ``trivial'' configurations of $\cF$, for example, those where every finite graph has an $\cF$-free orientation, and more precisely those appearing in item~\ref{label:tractCases1}.~and~\ref{label:tractCases2}.~of Corollary~\ref{cor:nonTechnicalVersion}---there is a ternary operation $f$, that allows to combine arbitrary $\cF$-free orientations $D_1, D_2, D_3$ of any given  graph $G$ into a new $\cF$-free orientation $f(D_1,D_2,D_3)$ of $G$. 
The operation $f$ is fully symmetric, that is, any permutation of its arguments yields the same output, and it additionally behaves like either a minority or a majority vote, i.e., either the equation $f(D_1,D_1,D_2)=D_2,$ or the equation $f(D_1,D_1,D_2)=D_1$ holds for all $\cF$-free orientations $D_1,D_2$ of an arbitrary graph $G$. We discuss these, so-called, \emph{minority} and \emph{majority operations} and their connection to the tractability of the CSP in detail in Section~\ref{sect:DecidableCriterion}.
In particular this ``symmetry'' of the solution space of $\CSP(\HF)$ implies that the latter is in P. In full  detail, using the technical notions found in the preliminaries, the following  dichotomy for the $\cF$-free orientation problem holds.
\color{black}

\begin{theorem}\label{sellingmainthm}
    Let $\cF$ be a finite set of forbidden tournaments. Then precisely one of the following holds:
    \begin{itemize}
        \item The polymorphism clone $\Pol(\HF)$ has a uniformly continuous minion homomorphism to $\sP$. In this case, $\CSP(\HF)$ is NP-complete by~\cite{wonderland}. 
        \item The polymorphism clone $\Pol(\HF)$  contains a ternary operation that is pseudo-cyclic modulo $\End(\HF)$. In this case, $\CSP(\HF)$ is in $P$.
    \end{itemize}
\end{theorem}

Our proof however yields a much more general result, thanks to the strength and flexibility of the new proof method.  Note that the edge relation of the undirected graph $\HF$ can be first-order defined in $\DF$ by the formula $(x\to y)\vee (y\to x)$, where $\to$ is the edge relation of the oriented graph $\DF$; we say that $\HF$ is a \emph{first-order reduct} of $\DF$. More generally, let $k\geq 2$ and consider a first-order definable $k$-ary relation $R$ on $\DF$ all of whose tuples induce cliques in $\HF$, a condition that is necessary for the crucial Lemma~\ref{lem:theNdominatedfunctionG}. 
Then using the homogeneity of $\DF$, it is not hard to see that $R$ is defined by a disjunction $\phi_1\vee \cdots\vee\phi_\ell$, where each $\phi_i$ is a conjunction describing a tournament on a $k$-tuple; in other words, $R$ consists of all $k$-tuples inducing in $\DF$ one of a finite number of allowed  tournaments. Let  $R_1,\ldots,R_m$ be relations of this form, and let $(\HF,R_1,\ldots,R_m)$ be the expansion of the graph $\HF$ by these relations. Then $\CSP(\HF,R_1,\ldots,R_m)$ is the problem where we are given an undirected graph $G$, together with constraints for each $i\in\{1,\ldots,m\}$ requiring that certain tuples of vertices of this graph must belong to $R_i$ (i.e., belong to the set of tournaments allowed by $R_i$); the question is whether the graph can be oriented to become an $\cF$-free oriented graph in such a way that all the constraints are satisfied.

\begin{theorem}\label{sellingmainthm_general}
    Let $\cF$ be a finite set of forbidden tournaments, and let $R_1,\ldots, R_m$ be relations first-order definable in $\DF$ consisting of tuples inducing tournaments.  Then precisely one of the following holds:
    \begin{itemize}
        \item The polymorphism clone $\Pol(\HF,R_1,\ldots, R_m)$ has a uniformly continuous minion   homomorphism to $\sP$. In this case, $\CSP(\HF,R_1,\ldots, R_m)$ is NP-complete by~\cite{wonderland}. 
        \item The polymorphism clone $\Pol(\HF,R_1,\ldots, R_m)$  contains a ternary operation that is pseudo-cyclic modulo $\End(\HF,R_1,\dots,R_m)$. In this case, $\CSP(\HF,R_1,\ldots, R_m)$ is in P.
    \end{itemize}
\end{theorem}

As an application, consider the \emph{$\cF$-free orientation completion problem}, where some edges of the input graph might already be oriented, and the task is to orient the remaining edges as to obtain an $\cF$-free oriented graph. Bodirsky and Guzm\'an-Pro in~\cite{forbiddenTournaments} showed that each such problem is equivalent to the CSP with a 2-element template, thus obtaining a P/NP-complete complexity dichotomy. We obtain this result as a direct consequence of Theorem~\ref{sellingmainthm_general} by considering the structure $(\HF,\mathrel{\to})$, where $\to$ denotes the binary relation defined by the formula $x \to y$.

\subsection{Related work}\label{sect:relatedwork}
Constraint satisfaction problems with a finite template are known to exhibit a P/NP-complete complexity dichotomy by the recent theorem of Bulatov and Zhuk~\cite{BulatovFVConjecture,Zhuk20}. If the template is infinite,  however, then there is provably no complexity dichotomy, even under the assumption of $\omega$-categoricity~\cite{BodirskyGrohe,GJKMP-conf}, a model-theoretic notion that can be viewed as an approximation of finiteness (see Section \ref{modeltheoryprelims} for a precise definition). It is however conjectured that the CSPs of certain infinite $\omega$-categorical structures, more specifically of \emph{first-order reducts of finitely bounded homogeneous structures}, do in fact exhibit a P/NP-complete complexity dichotomy.
For a modern formulation of the conjecture, originally posed by Bodirsky and Pinsker in 2012~\cite[Conjecture 1.2]{BPP-projective-homomorphisms}, we refer to~\cite[Conjecture 1.9]{wonderland}. Here we only provide the necessary information about the scope of the conjecture.
A homogeneous relational structure $\bB$ with finite relational signature $\tau$ is said to be \emph{finitely bounded} if there is some natural number $b_\bB$ such that a finite $\tau$-structure $\bF$ embeds into $\bB$ (i.e., is isomorphic to some substructure of $\bB$) if and only if every substructure $\bF' \subseteq \bF$ of size at most $b_\bB$ embeds into $\bB$. By taking $b_{\DF}\geq 2$ bigger than the size of any tournament in $\cF$, it follows that $\DF$ is a finitely bounded homogeneous structure, so the CSP of its first-order reduct $\HF$ falls into the scope of the Bodirsky-Pinsker conjecture. 

The conjecture has been verified for various classes of structures such as the first-order reducts of $(\bQ,<)$~\cite{tcsps-journal}, of any  homogeneous graph~\cite{BMPP16} including the random  graph~\cite{BodPin-Schaefer-both},  the universal homogeneous  poset~\cite{posetCSP18}, the universal homogeneous tournament~\cite{MottetPinskerSmooth}, certain homogeneous hypergraphs including the universal one~\cite{Hypergraphs}, and many more. The recent theory of smooth approximations~\cite{MottetPinskerSmoothConf,MottetPinskerSmooth}, pioneered by Mottet and Pinsker, provides a uniform framework for proving such complexity dichotomies which builds on elements that were  previously used sporadically; 
it has been used to obtain for the first time the last two results mentioned above, and to reprove the previous ones in a streamlined manner. Our complexity dichotomy constitutes progress on the Bodirsky-Pinsker conjecture for first-order reducts of finitely bounded homogeneous oriented graphs, while our proof  refines the smooth approximations approach along the way. The reason for the latter is that the automorphism group of $\DF$ (except in trivial cases) acts on pairs of distinct vertices  with three orbits (a pair $(x,y)$ can be non-adjacent, or  an arc, or an inverted arc). In contrast to this, the previous classifications above had at most two such orbits; on a technical level, this step is reminiscent of the challenge of passing from a two-element domain to a three-element domain in the realm of finite-domain CSPs~\cite{Bulatov-3-conf,Schaefer}.

At this point, we would like to point out that Bitter and Mottet independently (see~\cite{FellerPinsker} and the comment in~\cite{forbiddenTournaments}) also announced in the conference paper~\cite{AntoineZenoRivalPaper} a new proof of  the theorem of Bodirsky and Guzm\'{a}n-Pro, obtaining a similar statement as Theorem~\ref{sellingmainthm_general}.  They themselves apply the technique of smooth approximations to obtain a complexity dichotomy for certain reducts of the directed graph $\DF$. Their proof, perhaps demonstrating the universality of the technique of smooth approximations, appears to be  structurally similar, yet not identical, to the one presented here.
 
\subsection{Outline of the article and the proof}\label{sect:outline}
In Section~\ref{Prelims} we introduce the necessary notions from the universal-algebraic approach to infinite-domain CSPs and review a description of the automorphism group of $\HF$ from~\cite{ClassificationOfReductsOfDF}. 

We then turn to the proof of Theorem~\ref{sellingmainthm_general}, following the typical strategy using smooth approximations. Namely, in Section~\ref{sect:mccores} we first establish that the structure under consideration is either not \emph{novel}, roughly meaning that its CSP is the CSP of some simpler structure whose complexity has previously been classified (also using smooth approximations); or it is novel but can be assumed to be a \emph{model-complete  core}, i.e., in a sense itself the simplest template of its CSP. For the former class of structures, Theorem~\ref{sellingmainthm_general} is obtained immediately from previous classifications, and we continue with the latter, for which we strive to prove a stronger statement (Proposition~\ref{prop:MainThmForNovelStructures}). We then make several preparatory observations about the \emph{polymorphisms}, i.e., ``higher-ary symmetries'' of such structures in Section \ref{PrepSmoothApprox}. Finally, in Section~\ref{applyingsmoothapprox}, we show that one of two things happen: either the \emph{canonical} polymorphisms (a certain subset of all  symmetries) of the structure satisfy non-trivial \emph{identities}, in which case the tractability of the CSP follows from~\cite{Bodirsky-Mottet}; or they do not, in which the theory of smooth approximations is used to lift this property (locally) to \emph{all} polymorphisms, yielding hardness of the CSP by~\cite{Topo-Birk}.

In Section~\ref{sect:DecidableCriterion} we provide a decidable tractability criterion for the generalised orientation problems of Theorem~\ref{sellingmainthm_general} which does not use any  infinite-domain CSP jargon. The section also contains some examples separating our more general results from those of Bodirsky and Guzm\'{a}n-Pro. We finish with Section~\ref{sect:outlook} discussing future research directions.

\section{Preliminary facts}\label{Prelims}

All structures appearing in this article are tacitly assumed to be (at most) countable, and their signatures are assumed to be finite. 

\subsection{Model-theoretic notions}\label{modeltheoryprelims}
A relational structure $\bB$ is \textit{homogeneous} if every isomorphism between finite induced substructures of $\bB$ extends to an automorphism of $\bB$. 
If $\cF$ is a set of tournaments, then there exists a homogeneous oriented graph  whose (induced) subgraphs are, up to isomorphism,  precisely all $\cF$-free oriented graphs (see e.g.~\cite{Hodges}); this oriented graph is unique up to isomorphism and shall be denoted by $\DF$. It is \emph{$\omega$-categorical}, i.e., its automorphism group $\Aut(\DF)$ acts on its $k$-tuples (componentwise) with only finitely many orbits, for all $k\geq 1$. We denote the underlying set of vertices of $\DF$ by $V$ and its edge relation by $\to$; the flipped edge relation shall be denoted by $\xto$. 

The edge relation $E$ of the undirected graph $\HF$ is defined in $\DF$ by $\{(x,y)\;|\; (x\to y)\vee (x\xto y)\}$. This graph is not necessarily homogeneous, but $\omega$-categorical since $\Aut(\HF)\supseteq \Aut(\DF)$. In fact, any structure with a first-order definition in $\DF$ is $\omega$-categorical, for the same reason; such structures are called \emph{first-order reducts} of $\DF$. 
 We denote by $N\subseteq V^2$ the \emph{non-edges} of $\HF$, i.e., injective pairs  of vertices which are not adjacent in $\HF$. We say that two tuples $x,y$ of elements of $V$ of the same length  have the same \textit{edge-type} if the map sending the former onto the latter is a partial isomorphism on $\HF$. Only if  $\HF$ is homogeneous this is equivalent to $x,y$ belonging to the same $\Aut(\HF)$-orbit.

The class $\age(\DF)$, \emph{the age of $\DF$}, consisting of finite substructures of $\DF$ up to isomorphism, has \emph{free amalgamation}, i.e., it is closed under not necessarily disjoint unions. This implies (see e.g.~\cite{42})  that there exists a linear order $<$ on $V$ such that $(\DF,<):=(V;\to,<)$ is homogeneous and has precisely the class of linearly ordered $\cF$-free graphs as its substructures; it is also $\omega$-categorical. The order $<$ is isomorphic to the order of the rationals and the structure $(\DF,<)$ is unique up to isomorphism and called the \emph{free superposition} of $\DF$ with the order of the rationals. It enjoys the \emph{Ramsey property}~\cite{CSS-Ramsey}, a combinatorial property whose definition itself will not be needed; we shall make a single reference to a consequence of it in Section~\ref{subsect:topology}. 

Another consequence of $\age(\DF)$ being a free amalgamation class is that $\DF$ has \emph{no algebraicity} (see e.g.~\cite{Oligo}), i.e., for all finite $F\subseteq V$, every orbit of the pointwise stabilizer of $F$ within $\Aut(\DF)$, acting on $V\setminus F$, is infinite. 
As $\Aut(\HF) \supseteq \Aut(\DF)$, also $\HF$ has no algebraicity and more generally the same holds true for any first-order reduct of $\DF$.

\subsection{Universal-algebraic notions}\label{sect:universalAlgebra}
A \emph{polymorphism} $f$ of a relational structure $\mathbb A$ is a homomorphism from a finite power $\mathbb{A}^n$ of $\mathbb A$ to  $\mathbb A$; here $n$ is said to be the \emph{arity} of $f$.
Polymorphisms are closed under composition and every coordinate projection $\bA^n \to \bA$ for arbitrary $n$ is a polymorphism of $\bA$. This implies that the set $\Pol(\bA)$ consisting of all polymorphisms of $\bA$, called the \emph{polymorphism clone} of $\bA$, is a \emph{function clone}.

In general, a function clone $\sC$ on a fixed set $C$, its \emph{domain}, is a set  of finitary operations on $C$  containing all projections and which is closed under arbitrary composition. For $n\geq 1$, the set of $n$-ary functions of $\sC$ will be denoted by $\sC^{(n)}$. We also say that the clone $\sC$ \emph{acts on} $C$. But $\sC$ has various other actions: on any power $C^I$ of $C$, by componentwise evaluation, on any subset $S$ of $C$ that is invariant under the action of $\sC$, by restriction, and on $S/\mathord{\sim}$, for an equivalence relation $\sim$ on $S$ that is invariant under the action of $\sC$, by evaluating operations on representatives of the $\sim$-classes.  We write $\sC \curvearrowright C^I$, $\sC \curvearrowright S$ and $\sC \curvearrowright S/\mathord{\sim}$ for these actions, and call  $(S,\mathord{\sim})$ a \emph{subfactor} of $\sC$. 
A subfactor $(S,\mathord{\sim})$ of $\sC$ such that $S/\mathord{\sim}$ contains at least two elements and such that $\sC$ acts on $S/\mathord{\sim}$ by projections is called a \emph{naked set}. 

For $n\geq 1$, an $n$-ary operation $c$ is called \emph{cyclic} if it satisfies $c(a_1,\ldots,a_n)=c(a_2,\ldots,a_n,a_1)$ for all elements $a_1,\ldots,a_n$ of its domain; we indicate this by writing $c(x_1,\ldots,x_n)\approx c(x_2,\ldots,x_n,x_1)$. It is called \emph{pseudo-cyclic modulo $U$}, where $U$ is a set of unary functions on its domain, if there exist $u,v\in U$ such that $u\circ c(x_1,\ldots,x_n)\approx v\circ c(x_2,\ldots,x_n,x_1)$. 
A 6-ary operation $s$ is called \emph{Siggers} if $s(x,y,x,z,y,z)\approx s(y,x,z,x,z,y)$. \emph{Pseudo-Siggers operations modulo $U$} are defined accordingly. 
A ternary operation $m$ is called a \emph{majority} if $m(x,x,y)\approx m(x,y,x)\approx m(y,x,x)\approx x$ and a \emph{minority} if $m(x,y,y) \approx m(y,x,y) \approx m(y,y,x) \approx x$. An operation $f$ is \emph{idempotent} if $f(x,\ldots,x)\approx x$; a function clone is idempotent if all of its operations are.

A mapping $\xi \colon \sC \to \sD$ between two function clones is a \emph{clone homomorphism} if it preserves arities, maps the $i$th coordinate projection to the $i$th coordinate projection (of the same arity) and respects the composition of functions. Spelled out this means, for all $n$ and all $n$-ary functions $f\in \sC$ and all $g_1,\dots, g_n \in \sC$ of some arity $k$, it holds $\xi( f \circ (g_1,\dots,g_n) ) = \xi(f) \circ (\xi(g_1),\dots,\xi(g_n))$. Observe that for a subfactor $(S,\mathord{\sim})$ of $\sC$, the natural mapping $\xi \colon \sC \to (\sC \curvearrowright S/\mathord{\sim})$ is a clone homomorphism. A clone $\sC$ is said to be \emph{equationally trivial} if it homomorphically maps into the clone of projections on a two element domain, which we denote by $\sP$. Otherwise it is called \emph{equationally non-trivial}.

A mapping $\xi \colon \sC \to \sD$ between two function clones is a \emph{minion homomorphism} if it preserves arities and respects \emph{minors}, i.e.,
$\xi(f(x_{\sigma(1)},\ldots,x_{\sigma(n)}))= \xi(f)(x_{\sigma(1)},\ldots,x_{\sigma(n)})$ for all $f\in\sC$ and all functions $\sigma\colon \{1,\ldots,n\}\to \{1,\ldots,n\}$.

\subsection{Topology and canonicity}\label{subsect:topology}
For all $n\geq 1$, the set $V^{V^n}$ of all $n$-ary functions on $V$ is naturally  equipped with the \emph{topology of pointwise convergence}, where for $f\in V^{V^n}$ and $S\subseteq V^{V^n}$ we have that $f$ is contained in the closure $\overline{S}$ of $S$ if for all finite $F\subseteq V^n$ there exists $g_F\in S$ which agrees with $f$ on $F$. If $\bA$ is any structure on $V$, then $\Pol(\bA)^{(n)}$ is closed in $V^{V^n}$, and $\Aut(\bA)$ is closed in $\Sym(V)$, where the latter subspace of $V^V$ is equipped with the subspace topology. If $\bA$ is any $\omega$-categorical structure, then $V^{V^n}/\Aut(\bA)$, defined as the quotient space where $f,g\in V^{V^n}$ are equivalent if $f\in\overline{\{\alpha g(x_1,\ldots,x_n)\;|\; \alpha\in\Aut(\bA)\}}$, is compact (see e.g.~\cite{BP-canonical}).

A function $f\in V^{V^n}$ is \emph{canonical} with respect to $\bA$ if $f$ and $f(\beta_1(x_1),\ldots,\beta_n(x_n))$ are equivalent in above factor space for all $\beta_1,\ldots,\beta_n\in\Aut(\bA)$; this means that  $f$ preserves the equivalence relations given by the partitions of $k$-tuples into $\Aut(\bA)$-orbits, and hence $f$ acts on these orbits, for all $k\geq 1$. It is called \emph{diagonally canonical} if this is the case whenever  $\beta_1=\cdots=\beta_n$. Given two functions $f,g\in V^{V^n}$, we say that $g$ \emph{locally interpolates} $f$ modulo  $\Aut(\bA)$ if  $f\in\overline{\{\alpha g(\beta_1(x_1),\ldots,\beta_n(x_n))\;|\; \alpha,\beta_1,\ldots,\beta_n\in\Aut(\bA)\}}\; .$ \emph{Diagonal interpolation} is defined similarly, again by requiring $\beta_1=\cdots=\beta_n$. Since $(\DF,<)$ has the Ramsey property, every function $g\in V^{V^n}$   locally (diagonally) interpolates modulo $\Aut(\DF,<)$ a function which is (diagonally) canonical with respect to $\Aut(\DF,<)$~\cite{BP-canonical, BPT-decidability-of-definability}.

If $\sC$ is any function clone on $V$, then \emph{uniform continuity} of a clone (or minion)  homomorphism $\phi$ from $\sC$ into $\sP$ is understood with respect to the natural metric which induces the above topology; it is equivalent to the existence of a finite set $F\subseteq V$ such that $\phi(f)$ is already determined by the restriction of $f$ to the appropriate power of $F$, for any $f\in\sC$.

\subsection{Smooth approximations}
In \cite{MottetPinskerSmooth} different notions of \emph{approximations} are introduced. We will briefly recall the relevant ones here. 
Let $A$ be a set, $n\geq 1$, $S$ and $X$ subsets of $A^n$ with $S \subseteq X$ and ${\sim}$ some equivalence relation on $S$. An equivalence relation $\eta$ on $X$ \emph{approximates} ${\sim}$ if the restriction of $\eta$ to $X$ is a (possibly non-proper) refinement of~${\sim}$. Let $\sG$ be a permutation group acting on $A$ leaving the $\sim$-classes invariant as well as preserving $\eta$. We say that the approximation $\eta$ is 
\begin{itemize}
    \item \emph{presmooth (with respect to $\sG$)} if for every equivalence class $C$ of $\sim$ there is some equivalence class $C'$ of $\eta$ such that $C \cap C'$ contains two disjoint tuples (i.e., the two tuples do not have a common entry) in the same $\sG$-orbit;
    \item \emph{very smooth (with respect to $\sG$)} if $\sG$-orbit equivalence is a (possibly non-proper) refinement of $\eta$ restricted to $S$.
\end{itemize}

While the name might suggest otherwise, not every very smooth approximation is presmooth; a counterexample can easily be constructed by considering the trivial group acting on a set of size at least two. 
However, if one considers group actions by automorphism groups of $\omega$-categorical structures without algebraicity (which will always be the case in this paper, and the intended situation for the theory of smooth approximations), then  every very smooth approximation is indeed presmooth. We refer the reader to~\cite{MottetPinskerSmooth} for further details.

Very smooth approximations are of central importance to the step of lifting the property of ``satisfying only trivial identities'' from the clone of canonical polymorphisms $\sC$  to the clone of all polymorphism $\sD$ (locally), as mentioned in Section~\ref{sect:outline}. Roughly speaking, we will be in the situation that every function from $\sD$ locally interpolates some function of $\sC$, $(S,\sim)$ is a naked set of $\sC$, and $\eta$ is a very smooth approximation of $\sim$ invariant under $\sD$. Note that $\sim$ need not be invariant under $\sD$, and that therefore $\sD$ need not act on its classes or even its support (contrary to $\sC$). However, it does act on $\eta$-classes, and using  that the approximation is very smooth we can construct a uniformly continuous clone homomorphism from $\sD$ to  $\sC \curvearrowright(S,\sim)$ (``Fundamental theorem of smooth approximations''), implying that $\sD$ satisfies no more identities than $\sC$.  As presmooth approximations are only an intermediate step towards very smooth approximations we now exhibit some sufficient conditions under which such an upgrade can be performed. The following result is a variant of Lemma 8 in~\cite{MottetPinskerSmooth}.

\begin{lemma}\label{presmoothtoverysmooth}
    Let $\bF$ be a homogeneous relational structure whose age has free amalgamation and $\bA$ be a first-order reduct of $\bF$. Let $\sC\subseteq \Pol(\bA)$ be some clone containing $\Aut(\bF)$, preserving ${\not =}$ and $X\subseteq A^k$ for some $k \geq 1$, such that $\sC \curvearrowright X$ preserves $\Aut(\bF)$-orbit equivalence. If $(S,\mathord{\sim})$ is some naked set of $\sC\curvearrowright X$ with $\Aut(\bF)$-invariant $\sim$-classes and $\eta$ some $\sC$-invariant equivalence relation on $X$ which presmoothly approximates $\sim$ with respect to $\Aut(\bF)$, then there is a naked set $(S',\mathord{\sim}')$ of $\sC\curvearrowright X$ that is very smoothly approximated by $\eta$ with respect to $\Aut(\bF)$.
\end{lemma}
\begin{proof}
    Let $T$ be the set consisting of $a\in X$ such that there is disjoint $b\in X$ in the same $\Aut(\bF)$-orbit with $a \mathrel{\eta} b$. The set $T$ is $\sC$-invariant as ${\not =}$, the $\Aut(\bF)$-orbit equivalence on $X$ and $\eta$ are. Since the approximation $\eta$ is presmooth with respect to $\Aut(\bF)$, $T$ intersects every $\sim$-class non-trivially. Setting $S' := S \cap T$ and letting $\sim'$ be the restriction of $\sim$ to $S'$ it follows that $(S',\sim')$ is a naked set of $\sC$ with $\Aut(\bF)$-invariant equivalence classes. Moreover $\eta$ presmoothly approximates $\sim'$ with respect to $\Aut(\bF)$. 

    We are left showing that the approximation is very smooth with respect to $\Aut(\bF)$. To do so take $a\in S'$ and denote its $\Aut(\bF)$-orbit by $O$. Note that $\sC$ containing $\Aut(\bF)$ and acting on $X$ implies $O \subseteq X$.
    By construction of $S'$ there is $b\in O$, disjoint from $a$ with $a \mathrel{\eta} b$. Since $\age(\bF)$ has free amalgamation, there is $c \in O$ such that no relations hold between $a$ and $c$ and moreover $(a,b)$ and $(c,b)$ belong to the same $\Aut(\bF)$-orbit. The invariance of $\eta$ implies $c \mathrel{\eta} b$, hence $a \mathrel{\eta} c$. For arbitrary $s,t\in O$ we choose $r \in O$ such that no relations hold between $r$ and $s$ and between $r$ and $t$, which is again possible by free amalgamation. Homogeneity of $\bF$ implies that $(a,c)$, $(r,s)$ and $(r,t)$ all lie in the same $\Aut(\bF)$-orbit, whence $r \mathrel{\eta} s$ and $r \mathrel{\eta} t$, implying $s \mathrel{\eta} t$. In summary, $O$ is contained in a single $\eta$-class, and since $a$ was chosen arbitrarily the claim follows.
\end{proof}

\subsection{The automorphism group of \texorpdfstring{$\HF$}{the graph reduct of the Henson digraph}}
In~\cite{ClassificationOfReductsOfDF} Agarwal and Kompatscher give a description of the automorphism group of first-order reducts of $\DF$. This description turns out to be useful for our purposes (it is crucially used in Lemma \ref{lem:theNdominatedfunctionG} and Lemma \ref{DFcanimpliesHFcan}), so we detail the relevant results. 
For a set $F \subseteq \Sym(V)$, we  write $\langle F \rangle$ to denote the (topologically) closed subgroup of $\Sym(V)$ generated by $F$.

A \textit{flip of $\DF$} is any  function $fl\in \Sym(V)$ preserving $E$ and $N$ and flipping every arc  of $\DF$, i.e., every pair  in $\to$ is sent to a pair in  $\xto$. It is easy to see that whenever $fl$ is a flip of $\DF$, then  the group $\langle \Aut(\DF) \cup \{fl\} \rangle$ contains all  flips of $\DF$,  and that the inverse of any  flip is again a flip of $\DF$. 

A \textit{simple switch of $\DF$} is any $sw \in \Sym(V)$  preserving $E$ and $N$ such that there exists some $a \in V$ such that $sw$ inverts an arc of $\DF$ if and only if it is incident to $a$. Again, if $sw$ is some simple switch of $\DF$, then the closed subgroup $\langle \Aut(\DF) \cup \{sw\} \rangle$ contains all  simple switches of $\DF$, and the inverse of any simple switch is again a simple switch. As a generalisation of a simple switch we define for any finite, possibly empty set $A\subseteq V$ an \textit{$A$-switch of $\DF$} to be any  function $sw_A \in \Sym(V)$ that preserves $E$ and $N$ and that flips an arc of $\DF$ precisely if the arc is incident to a unique vertex in $A$. Note that $\emptyset$-switches are precisely the automorphisms of $\DF$.
If we do not specify $A$, we call any such function a \emph{switch}. 

\begin{theorem}\label{descriptionOfAutHF}\cite[Theorem 2.2]{ClassificationOfReductsOfDF}
    Let $\bA$ be a first-order reduct of $\DF$ with $\Aut(\bA) \subseteq \Aut(\HF)$ such that one of the following cases applies:
    \begin{itemize}
        \item the graph $\HF$ is not homogeneous;
        \item $\Aut(\bA)$ is a proper subgroup of $\Aut(\HF)$.
    \end{itemize}
    Then $\Aut(\bA)$ is among the following groups: $\Aut(\DF)$, $\langle \Aut(\DF) \cup \{fl\} \rangle$, $\langle \Aut(\DF) \cup \{sw\} \rangle$, $\langle \Aut(\DF) \cup \{fl, sw \} \rangle$.
\end{theorem}

\begin{corollary}\label{NormalFormAutomorphisms}
    Let $\bA$ be as in Theorem~\ref{descriptionOfAutHF}. Then for every $\alpha \in \Aut(\bA)$ and for every finite set $F\subseteq V$  there exists  
    a switch $\beta$
    and  a flip $\gamma$ such that $\alpha$ agrees with $\beta$ or with $\gamma\circ\beta$ on $F$.
\end{corollary}
\begin{proof}
    Note that every automorphism of $\DF$ is a $\emptyset$-switch. Hence, by Theorem~\ref{descriptionOfAutHF}, we have that $\alpha$ agrees with some permutation  $\beta_1\circ\cdots\circ \beta_n$ on $F$, where each $\beta_i$ is either a switch or a flip. Since the composition of two switches is a switch, and the composition of two flips is an automorphism of $\DF$, we may assume that there are no two consecutive switches or flips in $\beta_1\circ\cdots\circ \beta_n$. Moreover, the composition $\gamma\circ \delta$ of a switch $\gamma$ with a flip $\delta$ is equal to the composition $\delta'\circ\gamma'$ of a flip $\delta'$ with a switch $\gamma'$. The statement follows. 
\end{proof}

\section{Model-complete cores}\label{sect:mccores}

Henceforth, we let $\bA=(\HF,R_1,\ldots, R_m)$ be as in Theorem~\ref{sellingmainthm_general}. In a first step, we observe that we can assume $\bA$ to be a model-complete core unless its CSP falls into the scope of a known complexity classification. A structure $\bB$ is a \emph{model-complete core} if for every endomorphism $g$ of $\bB$ and every finite subset $F$ of its domain there is an automorphism of $\bB$ which agrees with $g$ on $F$. For every $\omega$-categorical structure $\bB$ there exists a model-complete core $\bC$ which has the same CSP as $\bB$; this structure is unique up to isomorphism~\cite{Cores-journal}. Any model-complete core $\bB$ has the property that all orbits of $\Aut(\bB)$ on tuples are invariant under $\Pol(\bB)$; in fact, this property is an equivalent definition of  model-complete core~\cite{Book}.

Recall that the following homogeneous structures are unique up to isomorphism, and that they are uniquely determined by their age: the random graph with age all finite graphs, for every $n>1$, the $n$th Henson graph with age all finite graphs that do not contain an $n$-clique, the countable clique with age all finite complete graphs, and the universal homogeneous tournament with age all finite tournaments.
\begin{lemma}\label{lem:cores}
    $\bA$ is a model-complete core,  or its model-complete core is either a one-element structure, or a first-order reduct of a homogeneous graph (either of the random graph, a Henson graph, or the countable clique), or of the universal homogeneous  tournament.
\end{lemma}
\begin{proof} 
    Let $\bC$ be the model-complete core of $\bA$. If $E=\emptyset$, then all relations of $\bA$ are empty, and $\bC$ is a one-element structure. Assume henceforth that $E\neq \emptyset$. By the results of~\cite{MottetPinskerCores}, summarized in Theorem~4 there, the Ramsey property of $(\DF,<)$ implies that $\bC$ can be assumed to be a first-order reduct of a homogeneous  substructure $\bB$ of $(\DF,<)$, and each relation of $\bC$ to be defined in $\bB$ by the same quantifier-free formula defining the corresponding relation of $\bA$ in $(\DF,<)$. In particular, $\bC$ is the structure which $\bA$ induces on the domain of $\bB$. Moreover, there exists an endomorphism $g$ of $\bA$ with the following properties:
\begin{enumerate}
    \item the range of $g$ is contained in $\bC$, and the age of the structure induced by the range of $g$ in $(\DF,<)$ equals $\age(\bB)$; 
    \item $g$ is $(\DF,<)$-canonical, and hence acts on the orbits of $\Aut(\DF,<)$;
    \item $g$ is \emph{range-rigid}, i.e., it acts as the identity function on those orbits of $\Aut(\DF,<)$ that intersect the range of $g$.
\end{enumerate}
    We briefly sketch how this follows from~\cite{MottetPinskerCores} since there it is only implicit. 
    The existence of an endomorphism $g$ with properties 1.~and 3.~follows immediately from the proof of~\cite[Theorem 4]{MottetPinskerCores}. Moreover, an inspection of the proof of~\cite[Lemma 15]{MottetPinskerCores} shows that we can additionally assume $g$ to be $(\DF,<)$-canonical. To be more precise, the range-rigid function whose existence is asserted in~\cite[Lemma 15]{MottetPinskerCores} is in fact constructed from a canonical function, and hence automatically $\mathscr G$-canonical (in our case $(\DF,<)$-canonical).
    
    We will first establish that $g$ is monotonic with respect to $<$. Let $O$ be the orbit of $\Aut(\DF,<)$ that is the image under $g$ of the $\Aut(\DF,<)$-orbit $U$ consisting of increasing (henceforth to be understood with respect to $<$) non-edges. Then $O$ consists of increasing pairs: if $W\subseteq E$ is any $\Aut(\DF,<)$-orbit in the range of $g$ consisting of increasing pairs, then there exist $x,y,z\in V$ with $(x,y),(y,z)\in W$ and $(x,z)\in U$; but then $(g(x),g(y)),(g(y),g(z))\in W$ by range-rigidity of $g$, so the pair $(g(x),g(z))\in O$ is increasing by transitivity.
    For arbitrary $x,z\in V$ with $x < z$, there is $y \in V$ with $(x,y),(y,z) \in U$, so $(g(x),g(y)), (g(y), g(z)) \in O$ imply $g(x) < g(z)$.  

    Assume first that the arcs in the range of $g$ are increasing only, or decreasing only. Let $\mathbb H$ be the graph induced by $\HF$ on the domain of $\bC$, and note that $\bH$ is a first-order reduct of $\bC$. The following argument shows that $\mathbb H$ is homogeneous and that $\Aut(\mathbb H)$ is contained in $\Aut(\bC)$, i.e., that $\bC$ is a first-order reduct of $\mathbb H$. 
    Let $u,v$ be isomorphic tuples in $\bH$. There are automorphisms $\alpha,\beta$ of $\DF$ such that $\alpha(u),\beta(v)$ are both increasing. We then have that $g(\alpha(u)),g(\beta(v))$ are isomorphic even in $(\DF,<)$, by the fact that $g$ is monotonic. Since $g\circ \alpha$ and $g\circ \beta$ act as endomorphisms on $\bC$, and $\bC$ is a model-complete core, there are automorphisms of $\bC$ which agree with the two functions on $u,v$, respectively. Since $g(\alpha(u)),g(\beta(v))$ are isomorphic in $(\DF,<)$, the homogeneity of $\bB$ then implies that they belong to the same orbit of $\Aut(\bB)$ and hence also of $\Aut(\bC)$. Putting this together, $u,v$ belong to the same orbit of $\Aut(\bC)$, so $\bH$ is a homogeneous graph as every automorphism of $\bC$ is an automorphism of $\bH$. Moreover the argument shows that $\Aut(\bH)$ is contained in $\Aut(\bC)$. 
    
    To determine the homogeneous graph $\bH$ we have to additionally distinguish whether $N$ intersects the range of $g$. Assume first that $N$ does not intersect the range of $g$, then $\age(\bH)$ consists of all finite complete graphs (use property (1) of $g$), and $\bH$ is the countable clique. 
    If $N$ intersects the range of $N$, the graph $\bH$ either is the random graph, or a Henson graph: applying $g$ to a tournament embedded in $\DF$ yields a transitive tournament as $g$ aligns the arcs with the order. So property (1) of $g$ implies that $\age(\bH)$ either consists of all finite graphs, if every transitive tournament is $\cF$-free; or it consists of all finite graphs that do not contain an $n$-clique, if the transitive tournament on $n$ vertices is the largest transitive $\cF$-free tournament. In the latter case an acyclic orientation of a graph is $\cF$-free if and only if the graph does not contain an $n+1$-clique and conversely no graph containing an $n+1$-clique has an $\cF$-free orientation.

    We now assume that the range of $g$ contains both increasing and decreasing arcs. Assume first that $N$ intersects the range of $g$. Then $(\DF,<)$ and $\bB$ are isomorphic by Fra\"isse's theorem: both are homogeneous relational structures and their age coincides, which readily follows from properties (1) and (3) of $g$. Hence $\bC$ is isomorphic to $\bA$, and $\bA$ is a model-complete core. Finally, assume that $N$ is not in the range of $g$.
    In this case $g$ maps a tuple in $N$ to a tuple in $\xto \cup \to$, without loss of generality sending increasing non-edges to increasing arcs, and the following argument implies that every tournament occurs in $\DF$. Suppose not every tournament occurs, then there is some smallest tournament $T$ with this property. The oriented graph $G$ one obtains by deleting one arc from $T$ is $\cF$-free by minimality of $T$. Now $(\DF,<)$ contains a copy $G'$ of $G$ such that the orientation of the deleted arc agrees with the order $<$, so $g$ maps the deleted arc to an arc of the same orientation. Consequently the underlying oriented graph on $g(G')$ is isomorphic to $T$, a contradiction.
    Now properties (1) and (3) of $g$, and an application of Fra\"iss\'e's theorem imply that $\bB$ is the universal homogeneous ordered tournament $(\mathbb T,<)$, so $\bC$ is a first-order reduct of $\mathbb T$, which is the universal homogeneous tournament, see e.g.~{\cite[Theorem 4.4.4]{Tent-Ziegler}}.
\end{proof}
 
We call $\bA$ \emph{novel} if its model-complete core is not a one-element structure, or a first-order reduct of either a homogeneous graph or of the universal homogeneous tournament; by the preceding lemma, this implies that $\bA$ is a model-complete core. This definition is inspired by the fact that if $\bA$ is not novel, then $\CSP(\bA)$ has already been classified in~\cite{BMPP16} and~\cite{MottetPinskerSmooth}, or is  trivial (if the core of $\bA$ is a one-element structure). We are going to prove Theorem~\ref{sellingmainthm_general} for novel $\bA$, which will be a tacit assumption from the next section on, and next argue how Theorem~\ref{sellingmainthm_general} follows from this case.

\begin{proposition}\label{prop:notnovel}
    Assume that $\bA$ is not novel. Then Theorem~\ref{sellingmainthm_general}  holds for $\bA$.
\end{proposition}
\begin{proof}
    Let $\bC$ be the model-complete core of $\bA$. It follows from~\cite{BMPP16} (in the case where $\bC$ is a first-order reduct of a homogeneous graph) and from~\cite{MottetPinskerSmooth} (in the case where $\bC$ is a first-order reduct of the universal homogeneous tournament) that either $\bC$ has a uniformly continuous clone homomorphism to $\sP$, or it has a ternary polymorphism which is pseudo-cyclic modulo $\overline{\Aut(\bC)}$ (for one-element $\bC$ this holds trivially). In the former case, $\bA$ has a uniformly continuous minion homomorphism to $\sP$  by~\cite{wonderland}, and hence the first case of Theorem~\ref{sellingmainthm_general} applies; in the latter case, the proof of~\cite[Corollary 6.2.]{Topo} shows that also $\bA$ has a pseudo-cyclic polymorphism modulo $\End(\bA)$.  
\end{proof}

We finish this section with two important observations that follow if $\bA$ is novel (so in particular a model-complete core).

\begin{lemma}\label{lem:novelAobservations}
If $\bA$ is a model-complete core, then $N$ and $\neq$ are invariant under $\Pol(\bA)$.
\end{lemma}
\begin{proof}
    $N$ is an orbit of both $\Aut(\DF)$ and  $\Aut(\HF)$, and $\Aut(\bA)$ lies between these two groups. Hence, $N$ is also an orbit of $\Aut(\bA)$. Since  $\bA$ is a model-complete core, it follows that $N$ is invariant under $\Pol(\bA)$.

    Since $E$ is preserved by $\Aut(\bA)$, there exists an orbit  $O\subseteq E$ of $\Aut(\bA)$. The formula $\exists z N(x,z)\wedge O(y,z)$ then holds for $(x,y)\in V^2$ if and only if $x\neq y$. Both $N$ and $O$ are preserved by $\Pol(\bA)$ since $\bA$ is a model-complete core. The above formula is a \emph{primitive positive definition} of $\neq$ from $N$ and $O$, and hence $\neq$ is preserved by $\Pol(\bA)$ as well, see~\cite{PPRussianGuys,PPGeiger}.
\end{proof}

\begin{lemma}\label{lem:unaryFunctionsNovel}
    If $\bA$ is a novel structure, then the range of every unary order-preserving function in $\Pol(\bA)$ that is canonical with respect to $(\DF,<)$ contains both decreasing and increasing edges.
\end{lemma}
\begin{proof}
    If $f$ was a counterexample to this claim, then applying the argument of the third paragraph of the  proof of Lemma~\ref{lem:cores} to $g\circ f$ instead of $g$ would contradict novelty.
\end{proof}

\section{Preparation to apply Smooth Approximations}\label{PrepSmoothApprox}

Recall that in this section the structure $\bA$ is assumed to be novel. Hence, $\Pol(\bA)$ preserves the binary relation $\not =$ by Lemma \ref{lem:novelAobservations}, and more generally, for $k \geq 1$, the set of injective $k$-tuples $I_k$.  

For compactness reasons, we denote the clone $\Pol(\bA)$ by $\sC$ and define the  following subclones of $\sC$:
\begin{itemize}
    \item $\sC_N$ is the subclone generated by the \textit{$N$-dominated} polymorphisms, where $f\in \sC^{(n)}$ is said to be $N$-dominated if for all $x_1,\dots,x_n \in V^2$ we have  $f(x_1,\dots,x_n) \in N$ whenever one of the $x_i$ is in $N$. More generally for any subclone $\sD\subseteq \sC$ we  write $\sD_N$ for the subclone of $\sD$ generated by the $N$-dominated functions in $\sD$; 
    \item $\CDF$ is the subclone of polymorphisms in $\sC$ which are canonical with respect to $\DF$;
    \item $\sC^{(\DF,<)}$ is the subclone of polymorphisms in $\sC$ which are canonical with respect to $(\DF,<)$;
    \item $\CDFinj$ is the subclone of polymorphisms in $\sC$ which preserve the equivalence of $\Aut(\DF)$-orbits of injective tuples; 
    \item $\sC^{(\DF,<),\mathrm{inj}}$ is the subclone of polymorphisms in $\sC$ which preserve the equivalence of $\Aut(\DF,<)$-orbits of injective tuples; 
    \item $\CAinj$ is the subclone of polymorphisms in $\sC$ which preserve the equivalence of $\Aut(\bA)$-orbits of injective tuples. 
\end{itemize}

Observe that $\CDF \subseteq \CDFinj$ since a function is contained in $\CDF$ if and only if it preserves $\Aut(\DF)$-orbit-equivalence of arbitrary tuples, while belonging to $\CDFinj$ is equivalent to preserving orbit-equivalence only for injective tuples. For the same reason we have
and $\CDF_N\subseteq \CDFinj_N$, $\sC^{(\DF,<)} \subseteq \sC^{(\DF,<),\mathrm{inj}}$ and $\sC^{(\DF,<)}_N \subseteq \sC^{(\DF,<),\mathrm{inj}}_N$. Apart from $\sD_N \subseteq \sD$, for any subclone $\sD\subseteq \sC$, there are no further obvious inclusions between the above listed clones in general. However, with the additional assumption of $\CDF$ being equationally trivial, we establish $\CDFinj_N \subseteq \CAinj_N$ and $\sC^{(\DF,<),\mathrm{inj}}_N \subseteq \CDFinj_N$ (see Lemma~\ref{DFcanimpliesHFcan} and~\ref{gettingMajority}), facts crucially needed in the proof of the main theorem.

The following lemma stipulates the existence of a special function contained in $\sC$, or rather for every $n\geq2$ a special $n$-ary function contained in $\sC$. Roughly speaking, this function allows to restrict ones attention from the set of all $\Aut(\DF)$-orbits of pairs $\{=,N,\xto,\to\}$, to those that are adjacent in $\HF$, $\{\xto,\to\}$ a property which is used in numerous places throughout this paper. The proof of Lemma~\ref{allThoseTrivialClones} describes a typical application.
The lemma does not require $\bA$ to be a novel structure.
\begin{lemma}\label{lem:theNdominatedfunctionG}
    For all $n\geq 2$ there is an $n$-ary  injective operation  $g\in\sC$ which
    \begin{itemize}
        \item is $(\DF,<)$-canonical,
        \item is monotonic as a function $(V^n,<_\mathrm{lex})\to(V,<)$, where $<_\mathrm{lex}$ denotes the lexicographic order induced by $<$,
        \item is $\DF$-canonical,
        \item acts as the  first projection on the two $\Aut(\DF)$-orbits $\{\xto,\to\}$ of adjacent pairs, and 
        \item is $N$-dominated, in fact it satisfies $g(X_1,\ldots,X_n)\subseteq N$ whenever $X_1,\ldots,X_n$ are $\Aut(\DF)$-orbits of pairs which are not all ${=}$ and not all contained in $E$. 
    \end{itemize}
    
    Any such  $g$ is contained in  $\CAinj_N$.
\end{lemma}
\begin{proof}
    We first construct a binary $g$ with the asserted properties (so $n=2$); the general case is then obtained by considering $g(x_1,g(x_2,g(\ldots,g(x_{n-1},x_n)))$.
    
    Let $D$ be the oriented graph  with vertex set $V^2$ and an arc from $(a,b)$ to $(x,y)$ if and only if there is an arc in $D_\cF$ from $a$ to $x$ and an edge in $H_\cF$ between $b$ and $y$. The oriented graph  $D$, equipped with the lexicographic order induced by $<$ on pairs, is an $\cF$-free linearly ordered oriented graph, so by a standard compactness argument there is an embedding $g$ of this structure into $(D_\cF,<)$.   One readily verifies that $g$ is a polymorphism of $\bA$; here we use that all relations of $\bA$ contain only tournaments. The five properties listed for $g$ follow immediately by construction. 
        
    We will now show that $g\in \CAinj_N$, i.e., $g$ preserves the $\Aut(\bA)$-orbit equivalence of injective tuples. Let $s,t,s',t'$ be injective tuples of equal length $k\geq 1$ with entries in $V$ and such that $s,s'$ and $t,t'$ lie in the same $\Aut(\bA)$-orbit, respectively.  We have to show that the (injective) tuples $r:=g(s,t)$ and $r':=g(s',t')$ belong to the same $\Aut(\bA)$-orbit. As $s,s'$ and $t,t'$ have the same edge-type respectively, so do $r,r'$, by the definition of $g$. If $\HF$ is homogeneous, then the novelty of $\bA$ implies that $\Aut(\bA)$ is a proper subgroup of $\Aut(\HF)$. Therefore the description of $\Aut(\bA)$ from Theorem~\ref{descriptionOfAutHF} applies. Note that if there is an arc between two coordinates of $r$, then $s$ has the same arc between the same coordinates; in other words, the arcs on $r$ are obtained by deleting some (possibly empty) set of arcs from $s$. The same holds true for $r'$ and $s'$. Let $\beta\in\Aut(\bA)$ be such that $\beta(s)=s'$. 
    By Corollary \ref{NormalFormAutomorphisms} there is some finite $A\subseteq V$ such that $\beta$ agrees with a switch $sw_A$, or such a switch composed with a flip, on $s$. Let $I \subseteq \{1,\dots,n\}$ be the set of indices such that $s_i \in A$ and define $B := \{ r_i \mid i \in I\}$. We then  have that $sw_B(r)$ or that tuple composed with a flip (if the flip was present in $\beta$), belongs to the same $\Aut(\DF)$-orbit as $r'$. Since $sw_B$ (or its composition with a  flip) belongs to $\Aut(\bA)$, this finishes the proof.   
\end{proof}

\begin{lemma}\label{nonTriviliatyNDominated}
    Let $\sD$ be a function clone whose functions preserve $N$.  If $\sD$ contains a pseudo-cyclic function modulo $\sD^{(1)}$, as well as an $N$-dominated function of the same arity, then   $\sD_N$ also contains a pseudo-cyclic function modulo $\sD^{(1)}_N$.
\end{lemma}
\begin{proof}
    Let $e,f,c\in\sD$ be so that $e\circ c(x_1,\ldots,x_n)\approx f\circ c(x_2,\ldots,x_n,x_1)$, and let $h(x_1,\ldots,x_n)$ be $N$-dominated. We have 
    \begin{align*}
        e &\circ c\big(h(x_1,\ldots,x_n),h(x_2,\dots,x_n,x_1),\ldots,h(x_n,x_1\ldots,x_{n-1})\big) \\
        \approx f &\circ c\big(h(x_2,\ldots,x_n,x_1),\ldots,h(x_n,x_1,\dots,x_{n-1}),h(x_1,\ldots,x_n)\big)
    \end{align*}
    and hence the $n$-ary function $c(h(x_1,\ldots,x_n),\ldots,h(x_n,x_1\ldots,x_{n-1}))$ is pseudo-cyclic modulo the unary part of $\sD_N$ as witnessed by $e,f\in \sD_N$. But this function  clearly is $N$-dominated as well.
\end{proof}

\begin{lemma}\label{allThoseTrivialClones}
    If one of the following clones $\CDF_N, \CDF, \CDFinj_N, \CDFinj$ is equationally trivial, then all of them are equationally trivial; this is the case if and only if the action of these clones on $\{\xto,\to\}$ is equationally trivial.
\end{lemma}
\begin{proof}
    Note that all mentioned clones  act on $\{\leftarrow, \rightarrow\}$, and these actions agree. Clearly, if this action is equationally trivial, then so are all mentioned clones. If it is  equationally non-trivial, then it contains a ternary cyclic function $c$~\cite[Theorem 4.2]{Cyclic}. Composing $c$ in its natural action with a ternary injection provided by Lemma~\ref{lem:theNdominatedfunctionG} in the same way as in the previous lemma yields an element of $\CDF_N$ which is pseudo-cyclic modulo  $\overline{\Aut(\DF)}$. Since $\CDF_N$ is contained in all other mentioned clones, all of them are equationally non-trivial.
\end{proof}

\begin{lemma}\label{DFcanimpliesHFcan}
    If $\CDF$ is equationally trivial, then $\CDFinj_N \subseteq \CAinj_N$.
\end{lemma}
\begin{proof}
    It suffices to prove that every $N$-dominated function $f\in \CDFinj$ is contained in $\CAinj$. Lemma \ref{allThoseTrivialClones} implies that the clone $\CDFinj\curvearrowright \{\xto,\to\}$ is equationally trivial, so it consists of essentially unary functions by Post's classification \cite{Post}. Hence for every $N$-dominated $f\in \CDFinj$ of arity $k\geq 1$, there exists $i \in \{1,\dots,k\}$ such that the action of $f$ on $\{ \xto,\to\}$ depends only on the $i$th coordinate. Note that for any injective tuples $a_1,\ldots,a_k$ of the same length all arcs of $f(a_1,\dots,a_k)$ are also found in $a_i$, just some arcs from $a_i$ might be missing.  
    Given tuples $b_1,\ldots,b_k$ from the same $\Aut(\bA)$-orbits as $a_1,\ldots,a_k$, we have that $f(b_1,\ldots,b_k)$ has the same edge-type as $f(a_1,\ldots,a_k)$. Let $\alpha \in \Aut(\bA)$ be an automorphism sending $a_i$ onto $b_i$. The same argument as in the proof of Lemma \ref{lem:theNdominatedfunctionG} shows that there is $\beta \in \Aut(\bA)$ such that $\beta \circ f(a_1,\dots,a_k)$ and $f(b_1,\dots,b_k)$ belong to the same $\Aut(\DF)$-orbit. This establishes that $f$ preserves the $\Aut(\bA)$-orbit equivalence of injective tuples.
\end{proof}

\section{Applying Smooth Approximations}\label{applyingsmoothapprox}

In this section we prove Theorem~\ref{sellingmainthm_general} for novel structures, which, together with Proposition~\ref{prop:notnovel}, establishes the theorem in full generality.
More precisely we will prove the following result. 

\begin{proposition}\label{prop:MainThmForNovelStructures}
    Let $\cF$ be a finite set of forbidden tournaments, and let $R_1,\ldots, R_m$ be relations first-order definable in $\DF$ consisting of tuples inducing tournaments such that the structure $\bA = (\HF,R_1,\ldots, R_m)$ is novel.  Then precisely one of the following holds:
    \begin{itemize}
        \item The polymorphism clone $\Pol(\bA)$ has a uniformly continuous clone   homomorphism to $\sP$. In this case, $\CSP(\bA)$ is NP-complete by~\cite{Topo-Birk}. 
        \item The polymorphism clone $\Pol(\bA)$  contains a ternary operation which is canonical with respect to $\DF$ as well as pseudo-cyclic modulo $\overline{\Aut(\DF)}$. In this case,  $\CSP(\bA)$ is in $P$, as follows from~\cite{Bodirsky-Mottet}.
    \end{itemize}
\end{proposition}

Note that in case of novel structures this result is indeed stronger than Theorem~\ref{sellingmainthm_general}. Namely every clone homomorphism is in particular a minion homomorphism, and, as $\Aut(\DF) \subseteq \End(\bA)$, every function that is pseudo-cyclic modulo $\overline{\Aut(\DF)}$ is also pseudo-cyclic modulo $\End(\bA)$. 

\subsection{The Master Proof of Proposition~\ref{prop:MainThmForNovelStructures}}\label{sect:masterproof}
If the clone $\CDF$ is equationally non-trivial, then so is the action of $\CDF_N$ on $\{\xto,\to\}$, by Lemma~\ref{allThoseTrivialClones}. Hence, $\CDF_N$ contains a ternary operation that acts as a cyclic operation on $\{\xto,\to\}$~\cite[Theorem 4.2]{Cyclic}. This operation is a pseudo-cyclic operation modulo $\overline{\Aut(\DF)}$ in its natural action, by a standard compactness argument, see e.g.~\cite[Lemma 3]{BP-canonical}.

So assume now that $\CDF$ is equationally trivial. By Lemma \ref{allThoseTrivialClones} this is equivalent to $\CDFinj_N$ being equationally trivial, and hence Lemma \ref{CHFIsTrivial} implies that $\CAinj_N$ is equationally trivial. Now Lemma \ref{trivalOnFiniteSet} provides a $k \geq 1$ such that the action $\CAinj_N \curvearrowright I_k/\Aut(\bA)$ is equationally trivial, where, by abuse of notation, $\Aut(\bA)$ denotes the orbit equivalence relation it induces on injective $k$-tuples. The assumptions of Theorem \ref{LoopLemma} are met, providing us with a naked set $(S,\mathord{\sim})$ of $\CAinj_N \curvearrowright I_k$ for which one of two cases applies.
        
In the first case $(S,\mathord{\sim})$ is approximated by a $\sC$-invariant equivalence relation $\eta$ on $I_k$ that is presmooth with respect to $\Aut(\bA)$. Observe that the approximation is also presmooth with respect to $\Aut(\DF)$: For each class $C$ of $\sim$ there is a class $C'$ of $\eta$ such that $C \cap C'$ contains two disjoint tuples $a,a'$ in the same $\Aut(\bA)$-orbit. Since $\DF$ has no algebraicity, there is $\alpha \in \Aut(\DF)$ fixing $a'$, such that $a$ and $b:=\alpha(a)$ are disjoint tuples (see e.g.~(2.16) of~\cite{Oligo}). We have $b\in C$ as the $\sim$-classes are in particular $\Aut(\DF)$-invariant and $b \in C'$ as the relation $\eta$ is in particular $\Aut(\DF)$-invariant. Hence the approximation is presmooth with respect to $\Aut(\DF)$.
Moreover $(S,\mathord{\sim})$ is a naked set of $\CDFinj_N\curvearrowright I_k$, as it is a subclone of $\CAinj_N \curvearrowright I_k$, by Lemma~\ref{DFcanimpliesHFcan}. An application of Lemma~\ref{presmoothtoverysmooth} yields a naked set $(S',\mathord{\sim'})$ of $\CDFinj_N \curvearrowright I_k$ with $\Aut(\DF)$-invariant $\sim'$-classes such that $\sim'$ is approximated by $\eta$ and such that the approximation is very smooth with respect to $\Aut(\DF)$. Proposition \ref{ConstructionOfCloneHom} then provides the desired uniformly continuous clone homomorphism from $\Pol(\bA)$ to the clone of projections $\sP$. 
        
The second case of Theorem~\ref{LoopLemma}, on the other hand, leads to a contradiction in the following way. Namely, it guarantees a \emph{weakly commutative function with respect to $(S,\mathord{\sim})$}, i.e., a binary function $f$ such that  $f(a,b) \sim f(b,a)$ for all $a,b \in I_k$ such that  $f(a,b),f(b,a)$ are in $S$ and disjoint.

Then Lemma \ref{lem:gettingMajFromCommutativeFunction} yields a majority operation contained in $\CDFinj \curvearrowright \{ \xto, \to \}$, so this clone is equationally non-trivial. By Lemma \ref{allThoseTrivialClones} this is equivalent to $\CDF$ being equationally non-trivial, contradicting our assumption of equational triviality of $\CDF$.

\subsection{Proofs of all claims of the Master Proof}
The following theorem combines two general theorems  from~\cite 
{MottetPinskerSmooth} and instantiates them for our purposes. 

\begin{theorem}\label{LoopLemma}
    Let $k\geq 1$, and suppose that $\CAinj_N \curvearrowright I_k/\Aut(\bA)$ is equationally trivial. Then there is a naked set $(S,\mathord{\sim})$ of $\CAinj_N \curvearrowright I_k$ with $\Aut(\bA)$-invariant $\sim$-classes such that one of the following holds:
    \begin{enumerate}
        \item The relation $\sim$ is approximated by a $\sC$-invariant equivalence relation $\eta$ on $I_k$ that is presmooth with respect to $\Aut(\bA)$.
        \item $\sC$ contains a \emph{weakly commutative function with respect to $(S,\mathord{\sim})$}.
    \end{enumerate}
\end{theorem}
\begin{proof}
    We apply Theorem~11 of \cite{MottetPinskerSmooth} to the clones  $\CAinj_N\subseteq \sC$ acting on $I_k$ (so in that theorem, $A=I_k$ and $n=1$). Then $\sC$ is the polymorphism clone of a model-complete core without algebraicity and the action of $\CAinj_N$ on $I_k$ is canonical with respect to the group $\Aut(\bA)$ of unary invertibles in $\sC$, by definition. We also have that the action of $\CAinj_N$ on $I_k/\Aut(\bA)$ is equationally trivial. The mentioned theorem thus indeed applies, and provides a naked set $(S,\mathord{\sim})$ of $\CAinj_N \curvearrowright I_k$ with $\Aut(\bA)$-invariant $\sim$-classes and two possible cases. 
    In the first case, it asserts that exactly item~1 above holds. In the second case, Theorem~11 of~\cite{MottetPinskerSmooth}  guarantees that every $\sC$-invariant binary symmetric relation $R$ on $I_k$ that contains a   pair $(a,b)\in S^2$ with $a \not \sim b$ and $a,b$ disjoint also  contains a pair $(c,c')$, where $c,c'$ belong to the same  $\Aut(\bA)$-orbit. Then item~2 above follows immediately from Lemma~14 of~\cite{MottetPinskerSmooth}.
\end{proof}

\begin{lemma}\label{gettingMajority}
    If $\sC^{(\DF,<)}_N \nsubseteq\CDFinj_N$, then the action $\CDFinj_N\curvearrowright \{\xto,\to\}$ contains a majority operation.
\end{lemma}
\begin{proof}
    Let $f$ be an $N$-dominated function contained in $\sC^{(\DF,<)}$ which is not contained in $\CDFinj$,  and denote its arity by $n$. As $f$ preserves the edge relation $E$ and the equivalence of $\Aut(\DF,<)$-orbits, it acts on the set $S = \{ ({\to} \cap {<}), ({\to} \cap {>}), ({\xto} \cap {<}), ({\xto} \cap {>}) \}$ of  $\Aut(\DF,<)$-orbits of adjacent  pairs. 

    Firstly we will demonstrate that we may assume that $f$ acts idempotently on $S$ while maintaining its distinctive properties. To do so, let $\hat f$ be the unary map defined by $x \mapsto f(x,\dots,x)$.
    By composing $f$ with a suitable self-embedding of $\DF$ from the outside, we may assume that $\hat f$ preserves the order $<$. Then $\hat f$ acts on $S':=\{({\to} \cap {<}), ({\xto} \cap {<})\}$, and its action on $S$ is determined by this action, e.g., if $\hat f$ maps $(\to \cap <)$ to $(\xto \cap <)$, then it maps $(\xto \cap >)$ to $(\to \cap >)$. 
    Suppose that the action of $\hat f$ on $S$ is not surjective; then $\hat f$ is constant on $S'$, which is impossible by Lemma~\ref{lem:unaryFunctionsNovel}. Therefore $\hat f$ acts as a permutation on $S'$, and as $S'$ consists of $2$ elements the function $\tilde f := \hat f \circ f$ acts idempotently on $S'$ and hence on $S$.
    Clearly, $\tilde f \in \sC^{(\DF,<)}$ is still $N$-dominated, and moreover still  not contained in $\CDFinj$, since $\hat f$ either flips all edges or preserves them, none of which alters the non-canonicity of $f$ when composed with $\hat f$. 
    Replacing $f$ with $\tilde f$ we may henceforth assume that $f$ acts idempotently on $S$. 

    Next we claim that for $i \in \{1,\dots,n\}$ there exist $o_i\in S'$ and $q_i\in S$ belonging to distinct $\Aut(\DF)$-orbits such that $f(o_1, \dots , o_n)$ and $f(q_1,\dots , q_n)$ belong to the same $\Aut(\DF)$-orbit. Indeed, otherwise we would have that for any $p_1, \dots , p_n \in S$, the $\Aut(\DF)$-orbit of 
    $f(p_1, \dots , p_n)$ would be determined by the value of $f(o_1, \dots , o_n)$, where $o_i$ is the unique element of $S'$
    which, viewed as a subset of $V^2$, is disjoint from the $\Aut(\DF)$-orbit generated by $p_i$, for all $ i \in \{1,\dots,n\}$. Hence, $f$ would act on the $\Aut(\DF)$-orbits in $E$, and being $N$-dominated, would therefore be an element of $\CDFinj$, a contradiction.

    Without loss of generality we may assume that $f(o_1,\dots,o_n) \subseteq \mathord{\to}$. By induction on $m\geq 2$ we now prove that for any $a,b,c \in E^m$ such that for every $i\in \{1,\dots,m\}$ at most one of the edges $a_i,b_i,c_i$ is in $\xto$, there is a ternary polymorphism $h$ of $\bA$ with $h(a,b,c) \in (\mathord{\to})^m$. 
        
    A standard compactness argument, which we now present for the convenience of the reader, then implies the existence of a ternary function $h \in \CDFinj$ acting like a majority operation on $\{\xto,\to\}$.
    Let $T$ be the set of finite subsets of $V^3$, partially ordered by $\subseteq$.
    For every $F\in T$ there is a ternary function $h^F \in \sC$ which acts like a majority operation on the edges within $F$. 
    Replacing $h^F$ by the function  $(x,y,z) \mapsto h^F\big(p^3_1(x,y,z), p^3_2(x,y,z), p^3_3(x,y,z)\big)$, where each $p^3_i$ is a ternary $N$-dominated function whose existence is stipulated in Lemma \ref{lem:theNdominatedfunctionG} with the $i$th coordinate taking the role of the first one, we may moreover assume that $h^F$ is $N$-dominated within $F$, and hence $h^F$ preserves the $\Aut(\DF)$-orbit equivalence of injective tuples within $F$.
    Now the net $(h^F)_{F \in T}$ yields the desired function $h$. More precisely, recall that $\Pol(\bA)^{(3)}/\Aut(\DF)$ is compact, so that the net $([h^F]_{\theta})_{F \in T}$ has an accumulation point and we show that any $h\in \sC$ contained in such an accumulation point preserves the equivalence of $\Aut(\DF)$-orbits of injective tuples and acts as a majority operation on $\{\xto,\to\}$. 
    Observe that $h$ has those desired properties, whenever it has the desired properties within every $F \in T$ and let $F \in T$  be such a finite set of triples.      
    As $[h]_\theta$ is an accumulation point of $([h^F]_\theta)_{F \in T}$, the definitions of the topology of pointwise convergence and of the quotient topology straightforwardly imply that there is $F' \in T$ with $F \subseteq F'$ and $\alpha \in \Aut(\DF)$ such that $\alpha \circ h^{F'}|_F = h|_F$. By construction $h^{F'}$ has the desired properties within $F'$, so in particular within $F$, and since the application of an automorphism of $\DF$ from the outside preserves the properties in question, also $h$ has the desired properties within $F$.

    We now return to the inductive proof of the statement made before the compactness argument.
    In the base case $m=2$ we can take $h$ to be a suitable projection. In the induction step let $m\geq 3$. By the induction hypothesis there are ternary $h,h',h'' \in \sC$ such that within $x := h(a,b,c)$ all but possibly the first coordinate $x_1$,  within $y := h'(a,b,c)$ all but possibly $y_2$,  and within $z := h''(a,b,c)$ all but possibly $z_3$  are elements of $\to$. As we are done otherwise, we may assume that $x_1$,  $y_2$, and $z_3$ indeed belong to $\xto$. Let $g$ be an  injective binary function as in  Lemma \ref{lem:theNdominatedfunctionG} acting as the the first projection on $\{\xto,\to\}$. This function allows us to assume that $x, y$ and $z$ have the same equality type, that is, two elements appearing in pairs in $x$ are equal precisely if the corresponding elements in $y$ and $z$ are equal: this is achieved by replacing $x$ by $g(x,g(y,z))$, $y$ by $g(y,g(x,z))$, and $z$ by $g(z,g(x,y))$. We may assume that there are no repetitions  in $x$ (nor in $y$ or $z$), i.e., $x_i\neq x_j$ for all distinct $i,j\in \{1,\dots,m\}$, as otherwise we are done by the induction hypothesis. Moreover, it is impossible that $x_i$ is obtained by flipping the entries of $x_j$: this is clear if $i,j\neq 1$, and if $i,j\neq 2$ by the same argument with $y$ (which is of the equality type as $x$), and if $i,j\neq 3$ by the same argument with $z$.
        
    This ensures that on $x$, viewed as subsets of $V$, there are two linear orders $\prec,\prec'$, which agree everywhere except on the entries of $x_1$, where they disagree. So the assumptions on the equality type of $x$ imply that $\prec$ and $\prec'$ in particular agree on $x_i$, for $i \in \{2,\dots,m\}$. Simply define a linear order on $x$ that makes the entries of $x_1$ the two largest elements. Then we can flip the order between the entries of $x_1$ without changing the order between the entries of $x$ that are not contained in $x_1$ to obtain $\prec$ and $\prec'$.
    By possibly flipping  both $\prec$ and $\prec'$ we may furthermore assume that the first entry of $x_2$ is smaller than the second entry of $x_2$ with respect to $\prec$ and $\prec'$. As $x$ and $y$ have the same equality type the same applies to $y$.  

    By the above discussion there are $u^k \in \{x,y\}$ as well as $\alpha_k \in \Aut(\DF)$ for all $k\in \{1,\dots,n\}$, such that $\alpha_k(u^k_1,u^k_2) \in q_k\times o_k$ and such that for a fixed $i\in \{3,\dots,m\}$ the tuples $\alpha_1(u^1_i),\dots, \alpha_n(u^n_i)$ belong to the same $\Aut(\DF,<)$-orbit. 
    More precisely choose $u^k \in \{x,y\}$ so that $u^k_1$ and $q_k$ have the same edge orientation. It follows that $u^k_2$ and $o_k$ have the same edge orientation as well. We then choose $\prec_k \in \{\prec,\prec'\}$ on $u^k$ so that the order it induces on $u^k_1$ equals the order of $q_k$. Moreover, by the previous paragraph the order $\prec_k$ is increasing on the tuple $u^k_2$ and hence agrees with the order of $o_k$. Since the edge and order relation of $(D_\cF,<)$ are freely superposed, there are $\alpha_k \in \Aut(\DF)$ such that $<$ restricts to $\prec_k$ on $\alpha_k(u^k)$ and the aforementioned properties of the $u^k$ and the $\alpha_k$ follow.
    
    Finally the choice of $o_1,\dots,o_n$ and $q_1,\dots,q_n$ as well as the fact that $f$ acts idempotently on $S$ imply that $f(\alpha_1 u^1,\dots,\alpha_n u^n)$ is an element of $(\to)^m$. Recalling that $u^k\in \{x,y\}$, and that $x$ and $y$ were obtained by applying some composition of the functions $h,h',h''$ as well as $g$ to $(a,b,c)$, we have constructed $f'\in \sC^{(3)}$ such that $f'(a,b,c) = f(\alpha_1 u^1,\dots,\alpha_n u^n) \in (\to)^m$, concluding the induction step.
\end{proof}

\begin{lemma}\label{CHFIsTrivial}
    If $\CDFinj_N$ is equationally trivial, then $\CAinj_N$ is equationally trivial as well.
\end{lemma}
\begin{proof}
    Suppose that $\CAinj_N$ is equationally non-trivial; we prove that $\CDFinj_N$ acts non-trivially on $\{\leftarrow,\rightarrow\}$. The lemma then follows by Lemma \ref{allThoseTrivialClones}.

    First observe that the action of $\CAinj_N$ on $\Aut(\bA)$-orbits of injective $k$-tuples is equationally non-trivial for every $k\geq 2$; moreover, it is  idempotent since $\bA$ is a model-complete core. Hence, for each $k\geq 2$ there is a function $s_k\in \CAinj_N$ satisfying the 6-ary Siggers identity in this action. Any  accumulation point $s$ of the sequence $(s_k)_{k\geq 2}$ in the 6-ary part of $\CAinj_N/\Aut(\bA)$ then satisfies the Siggers identity in all of these actions; such an accumulation point exists by compactness. In its standard action, the function $s$ then satisfies the pseudo-Siggers identity on injective tuples, i.e., for all injective tuples $t_1,t_2,t_3$ of equal length we have that $s(t_1,t_2,t_1,t_3,t_2,t_3)$ and $s(t_2,t_1,t_3,t_1,t_3,t_2)$ belong to the same $\Aut(\bA)$-orbit.

    As $(\DF,<)$ has the Ramsey property, we can apply Lemma~21 of \cite{BPT-decidability-of-definability} to derive   that there exists $s'\in\sC^{(\DF,<)}$ which is locally interpolated by $s$ modulo $\Aut(\DF,<)$. Further the function $s'$ is contained in $\sC^{(\DF,<)}_N$, as it is locally interpolated by the $N$-dominated function $s$. 
    Suppose that $s'$ is not contained in $\CDFinj$. In this case Lemma~\ref{gettingMajority} yields a majority operation contained in $\CDFinj \curvearrowright \{ \xto,\to\}$, which is equivalent to $\CDFinj_N$ acting non-trivially on $\{\xto,\to\}$ by Lemma~\ref{allThoseTrivialClones}. Henceforth we may assume that $s' \in \CDFinj$.
    We claim that $s'$ still satisfies the pseudo-Siggers identity on injective tuples modulo $\Aut(\bA)$, as above. Let injective tuples $t_1,t_2,t_3$ of equal length be given. On the finite set $S$ consisting of all entries of these tuples, we have 
    \begin{align*}
        s'(x_1,\ldots,x_6)\approx\beta s(\alpha_1 x_1,\ldots,\alpha_6 x_6)
    \end{align*}
    for suitable $\alpha_1,\ldots,\alpha_6,\beta\in\Aut(\DF,<)$. Moreover, since $s$ acts on $\Aut(\bA)$-orbits of injective tuples, we have 
    \begin{align*}
        \gamma s(t_1,t_2,t_1,t_3,t_2,t_3) &= s(\alpha_1 t_1,\alpha_2 t_2,\alpha_3 t_1,\alpha_4 t_3,\alpha _5 t_2,\alpha_6 t_3), \\
        \gamma'  s(t_2,t_1,t_3,t_1,t_3,t_2) &= s(\alpha_1 t_2,\alpha_2 t_1,\alpha_3 t_3,\alpha_4 t_1,\alpha_5 t_3,\alpha_6 t_2),
    \end{align*}
    for suitable $\gamma,\gamma'\in\Aut(\bA)$. Picking $\delta\in\Aut(\bA)$ such that 
    \begin{align*}
        s(t_1,t_2,t_1,t_3,t_2,t_3) = \delta s(t_2,t_1,t_3,t_1,t_3,t_2),
    \end{align*}
    we get 
    \begin{align*}
        s'(t_1,t_2,t_1,t_3,t_2,t_3) &= \beta s(\alpha_1 t_1, \alpha_2 t_2, \alpha_3 t_1, \alpha_4 t_3, \alpha_5 t_2, \alpha_6 t_3) \\
        &= \beta \gamma s(t_1,t_2,t_1,t_3,t_2,t_3) \\
        &= \beta \gamma \delta s(t_2,t_1,t_3,t_1,t_3,t_2) \\
        &= \beta \gamma \delta (\beta \gamma')^{-1} \beta \gamma' s(t_2,t_1,t_3,t_1,t_3,t_2) \\
        &= \beta \gamma \delta (\beta \gamma')^{-1} \beta s(\alpha_1 t_2,\alpha_2 t_1,\alpha_3 t_3,\alpha_4 t_1,\alpha_5 t_3,\alpha_6 t_2) \\
        &= \beta \gamma \delta (\beta \gamma')^{-1} s'(t_2,t_1,t_3,t_1,t_3,t_2),
    \end{align*}
    so $\beta \gamma \delta (\beta \gamma')^{-1} \in \Aut(\bA)$ proves our claim.

    Suppose that $\CDFinj_N$ was  equationally trivial. Then its action on $\{\leftarrow,\rightarrow\}$ would be equationally trivial as well by Lemma \ref{allThoseTrivialClones}, and hence it would be an action by essentially unary functions -- this follows from Post's classification of function clones on a 2-element domain \cite{Post}. In particular, this would be true for $s'$. The function $f(x):=s'(x,\ldots,x)$ then must be a permutation of $\{\leftarrow,\rightarrow\}$ in this action, by the assumed triviality of $\CDFinj_N$. Composing $s'$ with $f$ from the outside, we may then assume that $s'$ acts as a projection on $\{\leftarrow,\rightarrow\}$, without loss of generality the first one. Pick injective tuples $a,b$ of the same edge-type in different $\Aut(\bA)$-orbits. To see that this is possible, assume the contrary, namely, that every pair of injective tuples with the same edge-type belongs to the same $\Aut(\bA)$-orbit. This implies $\Aut(\HF) \subseteq \overline{\Aut(\bA)} \cap \Sym(V) = \Aut(\bA)$, so in fact $\Aut(\HF) = \Aut(\bA)$. Consequently $\HF$ is a homogeneous graph and $\bA$ is a first-order of the former, contradicting novelty of $\bA$.
    Then $s'(a,b,a,b,b,b)$ and $s'(b,a,b,a,b,b)$ belong to the same $\Aut(\bA)$-orbit by the previous paragraph. On the other hand, since $s'$ acts as a projection on $\{\leftarrow,\rightarrow\}$ and is $N$-dominated, $s'(a,b,a,b,b,b)$ belongs to the same $\Aut(\DF)$-orbit as $a$ and $s'(b,a,b,a,b,b)$ to the same $\Aut(\DF)$-orbit as $b$, a contradiction.
\end{proof}

The following lemma holds for every clone that is the polymorphism clone of a model-complete core, has a ternary injection, and acts on orbits of injective tuples.

\begin{lemma}\label{trivalOnFiniteSet}
    If $\CAinj_N$ is equationally trivial, then there is some $k\geq 2$ such that the action of $\CAinj_N$ on $\Aut(\bA)$-orbits of injective $k$-tuples is equationally trivial.
\end{lemma}
\begin{proof}
    Suppose that the action of $\CAinj_N$ on $\Aut(\bA)$-orbits of injective $k$-tuples is equationally non-trivial for every $k\geq 2$. Since each of this actions is idempotent, and the underlying set of the action is finite, it follows that for every $k\geq 2$ there is a function $s_k$ in $\CAinj_N$ which acts as a 6-ary Siggers operation on the $\Aut(\bA)$-orbits of injective $k$-tuples. Similarly as in Lemma~\ref{gettingMajority}, we consider the  factor space $\CAinj_N/\Aut(\bA)$, whose 6-ary part is compact. Thus, the sequence $([s_k])_{k\geq 2}$, where $[f]$ denotes the class of any  function $f$ in the factor, has an accumulation point $[s]$. We then have that $s$  acts as a 6-ary Siggers operation on $\Aut(\bA)$-orbits of injective $k$-tuples for all $k\geq 2$.

    Let $g\in \CAinj_N$ be a ternary injection provided by Lemma~\ref{lem:theNdominatedfunctionG} and set $u(x,y,z):=g(x,y,z)$, $v(x,y,z):=g(y,z,x)$, and $w(x,y,z):=g(z,x,y)$. We have that $s(u,v,u,w,v,w)(t_1,t_2,t_3)$ and $s(v,u,w,u,w,v)(t_1',t_2',t_3')$ belong to the same $\Aut(\bA)$-orbits whenever $t_1,t_2,t_3,t_1',t_2',t_3'$ are finite tuples of equal length such that $t_i,t_i'$ belong to the same $\Aut(\bA)$-orbit for all $1\leq i\leq 3$. This implies that there exist $e,f\in \overline{\Aut(\bA)}$ such that the identity 
    \begin{align*}
        e\circ s(u,v,u,w,v,w)(x,y,z)\approx f\circ s(v,u,w,u,w,v)(x,y,z) 
    \end{align*}
    holds in the standard action of $\CAinj_N$. This identity, when written using $g$ rather than $u,v,w$ (along with $e,f,s$), is  non-trivial, hence showing that $\CAinj_N$ is indeed  equationally non-trivial.  
\end{proof}

\begin{proposition}\label{ConstructionOfCloneHom}
    Let $k\geq 1$, $(S,\mathord{\sim})$ a naked set of $\CDFinj_N \curvearrowright I_k$ with $\Aut(\DF)$-invariant $\sim$-classes and $\eta$ a $\sC$-invariant equivalence relation on $I_k$ that very smoothly approximates $\sim$ with respect to $\Aut(\DF)$. Then there is a uniformly continuous clone homomorphism $\sC \to \sP$.
\end{proposition}
\begin{proof}
    It suffices to construct a uniformly continuous clone homomorphism $\sC\to \CDFinj_N \curvearrowright S/\mathord{\sim}$, since $(S,\mathord{\sim})$ is a naked set of $\CDFinj_N\curvearrowright I_k$. 

    First we will demonstrate that we can assume that all tuples in $S$ have the same edge-type. To do so, pick any $x,y\in S$ with $x\not\sim y$, and set $x':=g(x,y)$ and $y':=g(y,x)$, where  $g$ is as in Lemma \ref{lem:theNdominatedfunctionG}. Then $x',y'\in S$ have  equal edge-type and $x'\not\sim y'$ since $g$, being an element of $\CDFinj_N$, acts as a projection on  $\sim$-classes. Since $\sC$ preserves $E$ and $N$, the set of all $k$-tuples of the same edge-type as $x',y'$ is $\sC$-invariant; hence, restricting $S$ to such tuples yields a naked set with the desired property. 

    By Lemma~\ref{lem:theNdominatedfunctionG}, for every $n\geq 1$ and every $1\leq i \leq n$ there is an $n$-ary  $(\DF,<)$-canonical injection $p^n_i$ which follows the $i$th coordinate for the order $<$,  acts like the $i$th projection on $\{\leftarrow, \rightarrow\}$, and such that $p^n_i(a_1,\ldots,a_n)$ belongs to $N$ whenever $a_1\ldots,a_n$ are pairs not all belonging to $E$ and not all belonging to ${=}$. Note that in particular, $p^n_i$ then acts  as the $i$th projection on $S/\Aut(\DF)$.

    We now construct the clone homomorphism  $\phi \colon \sC \to \CDFinj_N \curvearrowright S/\mathord{\sim}$ as follows. Given an $n$-ary $f \in \sC$ we let $f'$ be a $(\DF,<)$-canonical function which is locally interpolated by $f$ modulo $\Aut(\DF,<)$. 
    The existence of such a function $f'$ follows from Corollary~6 of \cite{BP-canonical} together with the fact that $(\DF,<)$ is a Ramsey structure implying that its automorphism group is extremely amenable~\cite{Kechris}. 
    We then set $f'' = f'(p^n_1,\dots,p^n_n)$ which is an $N$-dominated $(\DF,<)$-canonical function. Moreover, $f''$ is injective: if $(a_1,\ldots,a_n)$ is a tuple of pairs in $V^2$ which do not all belong to $=$, then, independently of $i$, we have that  $p^n_i(a_1,\ldots,a_n)$ belongs to $E$ if and only if all of the $a_j$ do, and to $N$ otherwise. Hence, further application of $f$ yields either an element of $E$ or of $N$, since these two sets are invariant under $f$. We next claim that $f''\in \CDFinj_N$. Otherwise, Lemma \ref{gettingMajority} implies that $\CDFinj \curvearrowright \{\xto,\to\}$ contains a majority operation, so by Lemma \ref{allThoseTrivialClones}, $\CDF$ would be equationally non-trivial, contradicting the assumption of this proposition.
    Hence $f''$ acts on the set $S/\mathord{\sim}$.
    We now let $\phi(f)$ be the action of $f''$ on $S/\mathord{\sim}$, and claim that this assignment is well-defined. To see this, we will first argue that for any $u_1,\dots, u_n \in S$ we have 
    \begin{align*}
        f''(u_1,\dots,u_n)\, \sim \, f'(u_1,\dots,u_n) \, ({\sim}\circ \eta) \, f(u_1,\dots,u_n)\; .
    \end{align*}
    Indeed, for the first equivalence observe that since all elements of $S$ are of the same edge-type, we have that  $u_i$ and $p^n_i(u_1,\ldots,u_n)$ belong to the same $(\DF,<)$-orbit and hence to the same $\sim$-class for all $i$; since $f'$ is $(\DF,<)$-canonical, its application to the $u_i$ and the $p_i(u_1,\ldots,u_n)$ therefore  yields tuples in the same $(\DF,<)$-orbit, and hence in the same $\sim$-class; in particular we see that $f'(u_1,\ldots,u_n)\in S$.
    For the second equivalence, note that $f'(u_1,\ldots,u_n)=\alpha f(\beta_1 u_1,\ldots,\beta_n u_n)$ for some $\alpha,\beta_1,\ldots,\beta_n\in\Aut(\DF,<)$. Since the approximation $\eta$ is very smooth with respect to $\Aut(\DF)$, we have that $u_i$ and $\beta_i u_i$ are $\eta$-related for all $i$. 
    With $f$ preserving $\eta$, this implies that  $f(\beta_1 u_1,\ldots,\beta_n u_n)$ and $f(u_1,\ldots,u_n)$ are $\eta$-related, and since $\sim$-classes are $\Aut(\DF)$-invariant we moreover get that $f'(u_1,\dots,u_n)$ and $f(\beta_1 u_1,\dots,\beta_n u_n)$ are $\sim$-related, proving the the claim.

    It follows from the above that 
    \begin{align*}
        \phi(f)([u_1]_\sim,\dots,[u_n]_\sim) = [f''(u_1,\dots,u_n)]_\sim = [f(u_1,\dots,u_n)]_{\eta\circ\sim},
    \end{align*} 
    implying that $\phi(f)$ does not depend on the choice of $f'$.
    Observe moreover that since $S/\mathord{\sim}$ is finite,  $\phi$ is a uniformly continuous mapping; it remains to show that it is a clone homomorphism. Clearly,  $\phi$ preserves arities and projections,  so we are left showing that it preserves composition. Let $f\in \sC$ be $n$-ary, $g_1,\dots,g_n$ be $m$-ary and $c_1,\dots, c_m \in S/\mathord{\sim}$. We calculate setwise
    \begin{align*}
        \phi(f(g_1,\dots,g_n) )(c_1,\dots,c_m) &= \big[f \big(g_1(c_1,\dots,c_m),\dots,g_n(c_1,\dots,c_m) \big)\big]_{\eta\circ\mathord{\sim}} \\
        &= \big[f \big( [g_1(c_1,\dots,c_m)]_\eta,\dots,[g_n(c_1,\dots,c_m)]_\eta \big)\big]_{\eta\circ\sim} \\
        &= \big[f \big( [g_1(c_1,\dots,c_m)]_\eta \cap S ,\dots,[g_n(c_1,\dots,c_m)]_\eta \cap S \big)\big]_{\eta\circ\sim} \\
        &\subseteq \big[f \big( [g_1(c_1,\dots,c_m)]_{\eta \circ \sim} ,\dots,[g_n(c_1,\dots,c_m)]_{\eta\circ \sim} \big)\big]_{\eta\circ\sim} \\
        &= \phi(f)( \phi(g_1),\dots,\phi(g_n))(c_1,\dots,c_m),
    \end{align*}
    and since $\phi(f(g_1,\dots,g_n) )(c_1,\dots,c_m)$ and $\phi(f)( \phi(g_1),\dots,\phi(g_n))(c_1,\dots,c_m)$ are both equivalence classes of $\sim$ they have to be equal, proving that $\phi$ preserves composition. 
\end{proof}

Recall that given a subclone $\sD \subseteq \sC$ acting on $I_k$ and $(S,\sim)$ some subfactor of $\sD \curvearrowright I_k$, a binary function $f \in \sC$ is said to be weakly commutative with respect to $(S,\sim)$ if for all $x,y\in I_k$ such that the tuples $f(x,y),f(y,x)$ are disjoint and contained in $S$, it holds that $f(x,y) \sim f(y,x)$.

\begin{lemma}\label{lem:gettingMajFromCommutativeFunction}
    Suppose that $k\geq 1$ and that $(S,\mathord{\sim})$ is a naked set of $\CAinj_N\curvearrowright I_k$ with $\Aut(\bA)$-invariant classes. If 
    $\sC$ contains a binary function $f$ that is weakly commutative with respect to $(S,\mathord{\sim})$, then $\CDFinj \curvearrowright \{\xto,\to\}$ contains a majority operation.
\end{lemma}
\begin{proof}
    Let $g$ be a  binary polymorphism of $\bA$ as in  Lemma \ref{lem:theNdominatedfunctionG}. Observe that the binary polymorphism $h$ defined by  $(x,y) \mapsto f( g(x,y), g(y,x) )$ is  injective since $g$ is and since $\not =$ is preserved by all polymorphisms of $\bA$,  in particular by $f$. Moreover, $h$ is still weakly commutative: if for $a,b\in I_k$ we have $h(a,b),h(b,a)\in S$ and the tuples are disjoint, then setting $a':=g(a,b)$ and $b':=g(b,a)$ we obtain the disjoint tuples $f(a',b'),f(b',a')\in S$, and hence $h(a,b)=f(a',b')\sim f(b',a')=h(b,a)$.     
    Hence replacing $f$ by $h$  we can assume that $f$ is itself injective. As the distinctive property of $f$ is stable under diagonal interpolation, we can also assume that $f$ is diagonally canonical with respect to $\Aut(\DF,<)$.

    Let $O = (N \cap {<})$ and $O' = (N \cap {>})$. Diagonal canonicity of $f$ implies that within $O$ the function $f$  acts on the set $\{({\xto}\cap {<}), ({\xto}\cap {>}),({\to}\cap {<}),({\to}\cap {<})\}$ of $\Aut(\DF,<)$-orbits of adjacent  pairs, i.e., whenever $x,y \in E$ with $(x_i,y_j) \in O$ for $i,j\in \{1,2\}$, then the $\Aut(\DF,<)$-orbit of $f(x,y)$ only depends on the $\Aut(\DF,<)$-orbits of $x$ and $y$. The same holds true for $O'$. An application of Lemma \ref{gettingMajority} shows that we may assume that within $O$ and $O'$, $f$ also acts on the  set $\{\xto,\to\}$ of $\Aut(\DF)$-orbits of adjacent pairs:  
    otherwise, let $e_1,e_2$ be self-embeddings of $(\DF,<)$ such that  $(e_1(x),e_2(y)) \in O$ for all $x,y \in V$ (or in $O'$ for the other case). Then the function $f(e_1(x),e_2(y))$ is $N$-dominated and  
    contained in $\sC^{(\DF,<),\mathrm{inj}}$   but not in $\CDFinj$, so Lemma \ref{gettingMajority} implies that $\CDFinj\curvearrowright \{\xto,\to\}$ contains a majority and we are done.
    To see that such embeddings $e_1,e_2$ actually exist consider the ordered oriented graph $(D,<_D)$ consisting of two disjoint copies $(D_1,<_1)$ and $(D_2,<_2)$ of $(\DF,<)$ (i.e., there are no edges between $D_1$ and $D_2$) such that every vertex in $D_1$ is smaller than every vertex in $D_2$. Let $q_i\colon (\DF,<) \to (D_i,<_i)$ be isomorphisms for $i=1,2$ and let $e\colon (D,<_D) \to (\DF,<)$ be an embedding which exists by universality of $(\DF,<)$. Setting $e_i = e\circ q_i$ for $i=1,2$ then gives the desired embeddings.

    So we may henceforth  assume that within $O$ and $O'$, the function $f$ acts on $\{\xto,\to\}$. It is easy to see that this action is necessarily essentially unary in both cases as no other function on $\{\xto,\to\}$ can be induced by a function on $\DF$. Similarly, one sees that the two actions on $\{\xto,\to\}$ cannot be constant, and hence are either by a projection, or by a projection composed with the unique non-trivial permutation of $\{\xto,\to\}$. 
    Moreover, if one of the actions is not by a projection, then $\overline{\Aut(\bA)}$ contains a function which acts as the non-trivial permutation on $\{\xto,\to\}$,  obtained as $x \mapsto f(e_1(x),e_2(x))$, where $e_1$ and $e_2$ are self-embeddings as above. 
    In this case $\Aut(\bA)$-orbits are invariant under flipping all edges of the underlying oriented graph structure. 
     
    We will now demonstrate that the actions within $O$ and $O'$ on $\{\xto,\to\}$ cannot depend on the same argument. Striving for a contradiction, suppose without loss of generality that both actions depend on the first argument. Then for any $a,b \in S$ of the same edge-type and  with $(a_i,b_j) \in O$ for all $i,j \in \{1,\dots,k\}$, we have that $f(a,b)$ belongs to the $\Aut(\bA)$-orbit of $a$, and in particular is an element of $S$; similarly, since $(b_j, a_i) \in O'$ for all $i,j \in \{1,\dots,k\}$, we see that $f(b,a)$ belongs to the $\Aut(\bA)$-orbit of $b\in S$. This immediately contradicts $f(a,b)\sim f(b,a)$ if we additionally assume $a\not\sim b$. But such $a,b$ meeting our assumptions above exist, as they can be obtained by taking any $\sim$-inequivalent $a',b'\in S$, and setting 
    $a=e_1\circ g(a',b')$ and $b=e_2 \circ g(b',a')$, with $g,e_1,e_2$ as above; here we use that $g\in \CAinj_N$ acts as a projection on $\mathord{\sim}$-classes. Henceforth we may therefore  assume  that $f$ depends on the first argument on $O$ and on the second argument on $O'$.

    We now consider the \emph{atomic diagonal $(\HF,<)$-type} $B$ of two non-adjacent edges $x,y$, i.e., of $(x,y) \in E^2$ with $(x_i,y_j) \in N$ for $i,j\in \{1,2\}$. More precisely $B$ consists of all pairs of tuples $(u,v)$ such that $uv$ is in the range of a partial isomorphism of $(\HF,<)$ with domain $xy$.
    Diagonal canonicity of $f$ implies 
    that within each such type $B$, $f$ acts on the set $\{\xto,\to\}$. We can moreover assume that the action of $f$ on $\{\xto,\to\}$ within each such type is idempotent or constant, by replacing $f$ with $f(e_1\circ f(x,y), e_2 \circ f(x,y))$, where $e_1,e_2$ are self-embeddings of $\DF$ ensuring that the atomic diagonal $(\HF,<)$-type of $(a,b)$ is the same as the atomic diagonal $(\HF,<)$-type of $(e_1\circ f(a,b), e_2\circ f(a,b))$ for all pairs of non-adjacent edges $(a,b)$.
    For the convenience of the reader, we give the elementary  argument showing the existence of such $e_1$ and $e_2$. Let $D$ be the oriented graph obtained as the disjoint union of two copies  $D_1, D_2$ of $\DF$, i.e., there are no edges between $D_1$ and $D_2$. Let $s,t$ be sequences in $V$ such that $(s,t)$ is an enumeration of $V^2$.   We consider the sequence 
    $f(s,t)$ as a sequence $x$ in $D_1$ (via some isomorphism $q_1$ from $\DF$ onto  $D_1$), and also similarly as a sequence $y$  in $D_2$ (via some isomorphism $q_2$).  Now   $f$ being injective implies that there is a linear order $<_D$ on $D$ such that all inequalities  witnessed on the pair of sequences $(s,t)$ with respect to $<$ are also witnessed on $(x,y)$ with respect to $<_D$, i.e., for all natural numbers $i,j$ we have that $s_i < s_j$ implies $x_i <_D x_j$, that 
         $t_i < t_j$ implies $y_i <_D y_j$, that 
         $s_i < t_j$ implies $x_i <_D y_j$, and that 
         $t_i < s_j$ implies $y_i <_D x_j$. 
    The ordered oriented graph $(D,<_D)$ embeds into $(\DF,<)$ via an embedding $e$ and we can set $e_i=e\circ q_i$ for $i=1,2$. 
    
    Suppose that the  action of $f$ on $\{\xto,\to\}$ within the atomic diagonal $(\HF,<)$-type $B$ of a pair of non-adjacent edges was constant. Without loss of generality we may assume that within $B$ it holds $f(\to,\to) = f(\xto,\xto) = \mathord{\to}$.
    We claim that this implies that every injective tuple in $\DF$ is mapped by some unary function in $\sC$ to an injective tuple, all of whose edges point forward within the tuple.
    The same argument as in the proof of Lemma \ref{lem:unaryFunctionsNovel} then shows that $\bA$ is a first-order reduct of the homogeneous graph $\HF$ contradicting the assumption that $\bA$ is novel.
    We now establish the claim by induction on the length of $a$.  If the length of $a$ is $2$, we pick $\alpha,\beta \in \Aut(\DF)$ so that $(\alpha(a),\beta(a))$ is in $B$ implying that $x \mapsto f(\alpha(x),\beta(x))$ is the desired function. For the inductive step assume that the length of $a$ is $n+1$, for some $n\geq 2$. By the induction hypothesis we may assume that all edges within the first $n$ entries of $a$ point forward. Assume first  that for $(u,v)\in B$ it holds $u_1 < v_1$. Let $\alpha,\beta \in \Aut(\DF)$ so that for $a' = \alpha(a)$ and $a'' = \beta(a)$ it holds $(a'_i,a''_j) \in O$, for $1 \leq i,j \leq n$ and $((a'_i,a'_{n+1}),(a''_i,a''_{n+1})) \in B$, for $1\leq i \leq n$; this is possible since the imposed  conditions do not force a cycle on the tuple $a'a''$ with respect to the order $<$. Then, $x \mapsto f(\alpha(x),\beta(x))$ is the desired unary function in $\sC$. If for $(u,v)\in B$ it holds $u_1 > v_1$ the same argument with $O'$ instead of $O$ yields the desired conclusion. 
    
    In summary, we now know that within the atomic diagonal $(\HF,<)$-type of every pair of non-adjacent edges, $f$ acts on $\{\xto,\to\}$ as a projection or as a semilattice operation, setting the stage for the final step in the proof, the construction of a ternary function in $\CDFinj$ acting like a majority on $\{\xto,\to\}$.
    
    Just like in the proof of Lemma \ref{gettingMajority} it suffices to show by induction on $m\geq 2$ that for all tuples $a,b,c \in E^m$ such that for every $1\leq i \leq m$ at most one of the pairs $a_i,b_i,c_i$ is in $\xto$, there exists a ternary $h\in \sC$ such that $h(a,b,c) \in {(\to)^m}$. For the base case $m=2$ we can take $h$ to be a suitable ternary projection. For the induction step let $m\geq 3$ and $a,b,c \in E^m$. Clearly we can assume that there are no repetitions, i.e., $(a_i,b_i,c_i) \not = (a_j,b_j,c_j)$ for all $1\leq i < j \leq m$.
    Using the same trick as in the proof of Lemma \ref{gettingMajority} we may also assume that the kernels of $a,b$ and $c$ are identical. This implies that there cannot be $i,j \in \{1,\dots,m\}$ such that flipping the pair $a_i$ yields $a_j$, as otherwise the same would be true for $b_i$ and $b_j$ as well as for $c_i$ and $c_j$,  contradicting the fact that for all $1\leq i \leq m$ at most one of the pairs $a_i,b_i,c_i$ is in $\xto$.
    By the induction hypothesis we may assume that all but the second component of $a$ are in $\to$ and all but the first component of $b$ are in $\to$, while the assumptions about the kernels of $a$ and $b$ are still true. 
    By applying suitable automorphisms of $\DF$ to $a$ and $b$ respectively, we can moreover assume that no element of a pair in $a$ is adjacent to any element of a pair in $b$. 

    If no equalities hold between $a_1$ and $a_2$ then there is $\alpha \in \Aut(\DF)$ with $(\alpha(a_1),\alpha(b_1)) \in O$ and $(\alpha(a_2),\alpha(b_2)) \in O'$. 
    Moreover, $(\alpha(a_i), \alpha(b_i))$ is a pair of non-adjacent edges for $3\leq i \leq m$, and so $f$ acts idempotently on $\{\xto,\to\}$ within its order type.
    Hence we obtain $f(\alpha(a),\alpha(b)) \in (\to)^m$ proving the existence of the desired polymorphism.  
    
    If there are equalities between $a_1$ and $a_2$ we only demonstrate how to deal with the case in which  the second coordinate of $a_1$ equals the first coordinate of $a_2$ as the other cases may be treated similarly. Note that by our assumptions on the kernel of $a$ we know that the first coordinate of $a_1$ and the second coordinate of $a_2$ are distinct elements.
    In the following, call an  atomic diagonal $(\HF,<)$-type $B$ of a pair of non-adjacent edges which is \textit{crossing} if for $(u,v) \in B$ the order $<$ between $u_1$ and $v_1$ is opposite of the order $<$ between $u_2$ and $v_2$. 
    
    If $f$ acts as a semilattice operation on $\{\xto,\to\}$ within a crossing $B$, then by flipping all order relations that hold in $B$, obtain a crossing type where $f$ acting on $\{\xto,\to\}$ prefers $\to$ over $\xto$; hence without loss of generality we may assume this is the case  on $B$. First assume $B$ demands $u_1 < v_1$ for a pair $(u,v)$, and hence $u_2 > v_2$. Let $\alpha \in \Aut(\DF)$ so that $(\alpha(a_1),\alpha(b_1)) \in B$ and such that $(p,q)\in O'$ whenever $p$ is an entry of $a_2$ and $q$ is an entry of $b_2$. Then $f(\alpha(a),\alpha(b)) \in (\to)^m$ and we are done. If on the other hand $B$ demands $u_1 > v_1$, and hence $u_2<v_2$,  for a pair $(u,v)$, then we pick  $\alpha \in \Aut(\DF)$ so that $(\alpha(b_1),\alpha(a_1)) \in B$ and $(\alpha(b_2),\alpha(a_2) )\in O$. Again,  $f(\alpha(b),\alpha(a)) \in (\to)^m$, finishing the induction step.    

    It remains to consider the case in which $f$ acts as a projection on $\{\xto,\to\}$ within each crossing type of non-adjacent edges. If $f$ behaves like the first projection within some crossing type that demands $u_1<v_1$, and hence $u_2>v_2$,  for a pair $(u,v)$, we may choose $\alpha \in \Aut(\DF)$ so that $(\alpha(a_1),\alpha(b_1)) \in B$ and $(\alpha(a_2),\alpha(b_2)) \in O'$ to again conclude  $f(\alpha(a),\alpha(b)) \in {(\to)^m}$. 
    Otherwise $f$ behaves like the second projection on $\{\xto,\to\}$ within each crossing type demanding $u_1<v_1$ for a pair $(u,v)$. Take such a crossing type $B$ and let $\alpha \in \Aut(\DF)$ such that $(\alpha(a_1),\alpha(b_1) ) \in O$ and $(\alpha(a_2) ,\alpha(b_2)) \in B$. Again $f(\alpha(a),\alpha(b)) \in (\to)^m$. 
\end{proof}

\section{Examples and a decidable tractability criterion}\label{sect:DecidableCriterion}

The goal of this section is to provide a tangible and in fact algorithmically decidable tractability criterion for the generalised orientation problems of Theorem~\ref{sellingmainthm_general} which moreover avoids all infinite-domain CSP jargon. Throughout the section, we let $\cF$ be a finite set of forbidden tournaments, and denote by $m_\cF$ the maximum number of vertices of a tournament in $\cF$. 
Let $R_1,\dots,R_m$ be relations first-order definable in $\DF$ consisting of injective tuples inducing tournaments in $\DF$ and let $n_i$ be the arity of $R_i$ for each $i$. We set $\bA = (\HF,R_1,\dots,R_m)$ and since we are interested in the computational complexity of its CSP we may, without loss of generality, assume that each $R_i$ is non-empty.

\begin{definition}
    Let $R$ be an $n$-ary relation on the domain of $\DF$ and invariant under $\Aut(\DF)$. The \emph{standard representation} of $R$ consists of all oriented graphs on $\{1,\ldots,n\}$ obtained by pulling back the relations induced by a tuple $x\in R$ in $\DF$ via the mapping $\{1,\ldots,n\}\to \DF$, $i\mapsto x_i$.
\end{definition}
Note that $\bA$ is uniquely determined by $\cF$ and the standard representation of each $R_i$, and thus has a finite representation from which we will decide tractability of the CSP.

\subsection{Distinguishing the old from the new}
Since we aim to determine tractability even for non-novel structures, we first need to establish how to recognize these. 
Denote the model-complete core of $\bA$ by $\bA'$ and recall that both structures have the same $\CSP$. 
We denote the countable clique by $K_\omega$, for $n > 2$ the $n$th Henson graph by $H_n$, and the universal homogeneous tournament by $\bT$.
\begin{lemma}\label{lem:nonNovelExplicit}
    The structure $\bA$ is novel unless one of the following cases holds.
    \begin{enumerate}
        \item \label{label:nonNovel1} The tournament on two vertices is contained in $\cF$. Then $\bA'$ is a one-element structure.      
        \item \label{label:nonNovel2} Every transitive tournament is $\cF$-free and the standard representation of every $R_i$ contains all transitive tournaments on $\{1\dots,n_i\}$.
        Then $\bA' = (K_\omega,R_1',\dots,R_m')$, where each $R_i'$ consists of all injective $n_i$-tuples of $K_\omega$.
        
        \item \label{label:nonNovel3} There exist $n>2$ such that every transitive tournament of size $<n$ is $\cF$-free and no tournament of size $\geq n$ is $\cF$-free. Moreover the standard representation of every $R_i$ contains all transitive tournaments on $\{1\dots,n_i\}$.
        Then $\bA' = (H_n,R_1',\dots,R_m')$, where each $R_i'$ consists of all injective $n_i$-tuples inducing cliques in $H_n$.
        
        \item \label{label:nonNovel4} Every tournament is $\cF$-free. Then $\bA' = (K_\omega,R_1',\dots,R_m')$, is a first-order reduct of $\bT$ and each $R_i'$ is defined in $\bT$ by the same quantifier-free formula as $R_i$. 
    \end{enumerate}
    In cases~\ref{label:nonNovel1}.-\ref{label:nonNovel3}., $\CSP(\bA)$ is in $P$.
\end{lemma}
\begin{proof}
    The statements straightforwardly follow by recalling the proof of Lemma~\ref{lem:cores}. Item~\ref{label:nonNovel1}.~corresponds to the situation $E=\emptyset$. 
    
    Item~\ref{label:nonNovel2}.~corresponds to the situation where the range of $g$ contains increasing only or decreasing only arcs, and moreover does not intersect $N$. The following argument shows that for all $i\leq m$ any $x \in V^{n_i}$ that induces a transitive tournament in $\DF$ is contained in $R_i$. Without loss of generality we may assume that the range of $g$ contains only increasing arcs of $(\DF,<)$. Pick any $y\in R_i$. By definition of the structure $(\DF,<)$ there is $\alpha \in \Aut(\DF)$ such that $<$ orders the entries of the tuple $\alpha(y)$ in the same way as $\to$, the edge relation of $\DF$, orders the entries of $x$. As $g$ aligns the edge relation $\to$ with $<$, $x$ and $g\circ \alpha(y)$ belong to the same $\Aut(\DF)$-orbit and the claim follows since $R_i$ is invariant under $g$ as well as under $\Aut(\DF)$.  
    
    Item~\ref{label:nonNovel3}.~corresponds to the situation where the range of $g$ contains increasing only or decreasing only arcs and additionally intersects $N$. The homogeneous graph $\bH$ cannot be the random graph in this case as this would contradict the fact that $\bA'$ is a model-complete core. The statement about the $R_i$ containing every tuple inducing a transitive tournament in $\DF$ follows exactly as above. 

    Finally, item~\ref{label:nonNovel4}.~corresponds to the situation where the range of $g$ contains both increasing and decreasing arcs and does not intersect $N$. 
\end{proof}

\subsection{Minority and majority operations}\label{sect:minorityMajority}
We now develop a suitable setup for the   tractability test for novel structures.  
Let $G=(V_G,E_G)$ be a graph, and let $D_1=(V_G,\to_1)$, $D_2=(V_G,\to_2)$, and $D_3=(V_G,\to_3)$ be orientations of $G$. 
The ternary \emph{minority} and \emph{majority operation} that combine the oriented graphs $D_1,D_2,D_3$ into a new orientation of $G$ are defined as follows. 
Both, $\mino(D_1,D_2,D_3)$ and $\majo(D_1,D_2,D_3)$ have vertex set $V_G$. The antisymmetric edge relation
of the oriented graph $\mino(D_1,D_2,D_3)$ contains the tuple 
$(u,v)\in V_G^2$ precisely when $(u,v)$ is contained in exactly one or three of the relations $\to_i$, for $i=1,2,3$. The antisymmetric edge relation of $\majo(D_1,D_2,D_3)$ contains the tuple 
$(u,v)\in V_G^2$ precisely when $(u,v)$ is contained in exactly two or three of the relations $\to_i$, for $i=1,2,3$, see Figure~\ref{fig:majoClosedRelation} for an example. Note that $\mino$ ($\majo$) acts as a minority (majority) operation, as defined in Section~\ref{sect:universalAlgebra}, on the 
set of all orientations of a given graph. Moreover, it is totally symmetric in its arguments, i.e., any permutation of the arguments yields the same output.

We say that a set $\cD$ of orientations of some graphs $G_1,\dots,G_n$ is \emph{preserved by the minority (majority) operation} whenever applying $\mino$ ($\majo$) to  elements of $\cD$ that orient the same graph $G_i$ yields an oriented graph contained in $\cD$. By the above equations this condition only needs to be checked for all triples of \emph{distinct} elements of $\cD$. 
We say that a relation  $R$ on $\DF$ that is invariant under $\Aut(\DF)$ is preserved by the minority (majority) operation if its standard representation  is preserved by the minority (majority) operation.  

Recently, Bodirsky and Guzm\'{a}n-Pro showed that for 
every $n\leq m_\cF$, every set of $\cF$-free tournaments on $\{1,\dots,n\}$ that is invariant under the majority operation is either empty or contains every tournament, see~\cite[Lemma 25]{forbiddenTournaments}. In particular any set of $\cF$-free tournaments that is preserved by the majority operation is also preserved by the minority operation. However, this does not imply that any \emph{relation} that is preserved by the majority operation is also preserved by the minority operation. A simple counterexample can be constructed as follows. 

Let $\cF$ be any set of forbidden tournaments such that both tournaments on $3$ vertices are $\cF$-free. Consider the ternary relation $R$ of $\DF$ defined by the formula $\phi_1 \vee \phi_2 \vee \phi_3$, where
\begin{align*}
    \phi_1 &:=(x_1 \to x_2) \land (x_3 \to x_1) \land (x_2 \to x_3), \\
    \phi_2 &:=(x_1 \to x_2) \land (x_3 \to x_1) \land (x_2 \to x_3),  \\
    \phi_3 &:= (x_1 \to x_2) \land (x_3 \to x_1) \land (x_2 \to x_3). 
\end{align*} 
Letting $R_i$ be defined by $\phi_i$, the unique tournament contained in the standard representation of each $R_i$ is depicted in Figure~\ref{fig:majoClosedRelation}. The figure moreover demonstrates that $R$ is preserved by the majority operation, while it is not preserved by the minority operation. 

\begin{figure}
\centering
\begin{tikzpicture}

  \begin{scope}[xshift = -6cm]
    \node [vertex] (1) at (0,0) {};
    \node at (1) [below left] {1};
    \node [vertex] (2) at (1,0) {};
    \node at (2) [below right] {2};
    \node [vertex] (3) at (0.5,0.86602540378) {};
    \node at (3) [above] {3};
    \node (L1) at (0.5,-0.7) {$a \in R_1$};
      
    \foreach \from/\to in {1/2, 3/1, 2/3}     
    \draw [arc] (\from) to (\to);
  \end{scope}

  \begin{scope}[xshift = -3.5cm]
    \node [vertex] (1) at (0,0) {};
    \node at (1) [below left] {1};
    \node [vertex] (2) at (1,0) {};
    \node at (2) [below right] {2};
    \node [vertex] (3) at (0.5,0.86602540378) {};
    \node at (3) [above] {3};
    \node (L1) at (0.5,-0.7) {$b \in R_2$};
      
    \foreach \from/\to in {2/1, 1/3, 2/3}     
    \draw [arc] (\from) to (\to);
  \end{scope}

  \begin{scope}[xshift = -1.cm]
    \node [vertex] (1) at (0,0) {};
    \node at (1) [below left] {1};
    \node [vertex] (2) at (1,0) {};
    \node at (2) [below right] {2};
    \node [vertex] (3) at (0.5,0.86602540378) {};
    \node at (3) [above] {3};
    \node (L1) at (0.5,-0.7) {$c \in R_3$};
      
    \foreach \from/\to in {2/1, 3/1, 2/3}     
    \draw [arc] (\from) to (\to);
  \end{scope}

  \begin{scope}[xshift = 2.5cm]
    \node [vertex] (1) at (0,0) {};
    \node at (1) [below left] {1};
    \node [vertex] (2) at (1,0) {};
    \node at (2) [below right] {2};
    \node [vertex] (3) at (0.5,0.86602540378) {};
    \node at (3) [above] {3};
    \node (L1) at (0.5,-0.7) {$\mino(a,b,c) \not \in R$};
      
    \foreach \from/\to in {1/2, 1/3, 2/3}     
    \draw [arc] (\from) to (\to);
  \end{scope}

  \begin{scope}[xshift = 6cm]
    \node [vertex] (1) at (0,0) {};
    \node at (1) [below left] {1};
    \node [vertex] (2) at (1,0) {};
    \node at (2) [below right] {2};
    \node [vertex] (3) at (0.5,0.86602540378) {};
    \node at (3) [above] {3};
    \node (L1) at (0.5,-0.7) {$\majo(a,b,c) \in R$};
      
    \foreach \from/\to in {2/1, 3/1, 2/3}     
    \draw [arc] (\from) to (\to);
  \end{scope}
  
\end{tikzpicture}
\caption{The relation $R$ is preserved by $\majo$ but not by $\mino$.}
\label{fig:majoClosedRelation}
\end{figure}

We finish this subsection by a lemma which connects preservation by a minority or majority with the existence of certain canonical functions, i.e., infinite objects. This will allow us to provide a finitary tractability criterion in the case of novel structures via Proposition~\ref{prop:MainThmForNovelStructures} later.

\begin{lemma}\label{lem:buildingCanonicalFunction}
    The following statements are equivalent. 
    \begin{itemize}
        \item There is a ternary function $f\in \Pol(\bA)$ that is $\DF$-canonical and acts a a minority (majority) operation on the set of $\DF$-orbits of edges in $\HF$, $\{\xto,\to\}$. 
        \item The set of $\cF$-free tournaments of size $\leq m_\cF$ as well as the relations $R_1,\dots,R_m$ are preserved by the minority (majority) operation. 
    \end{itemize}
\end{lemma}
\begin{proof}
    The second statement follows from the first statement by unravelling all definitions. Conversely, the first statement follows from the second by constructing, similarly to the construction of $g$ in Lemma~\ref{lem:theNdominatedfunctionG}, a ternary $\DF$-canonical function $f$ that acts as a minority (majority) on $\{\xto,\to\}$, and additionally satisfies $f(X_1,X_2,X_3)\subseteq N$ whenever $X_1,X_2,X_3$ are $\Aut(\DF)$-orbits of pairs which are not all ${=}$ and not all contained in $E$. Note that the assumption of all $\cF$-free tournaments of size $\leq m_\cF$ being preserved by minority (majority) is in fact equivalent to the statement that the set of all $\cF$-free tournaments is preserved by minority (majority).  
\end{proof}

\subsection{The decidable tractability criterion}

We are now ready to provide a decidable reformulation of the tractability criterion from Theorem~\ref{sellingmainthm_general} avoiding all infinite-domain CSP jargon.

\begin{corollary}\label{cor:nonTechnicalVersion}
Either $\CSP(\bA)$ is NP-complete, or one of the following cases holds, in all which $\CSP(\bA)$ is in P. 
    \begin{enumerate}
        \item \label{label:tractCases1} Every transitive tournament is $\cF$-free and the standard representation of every $R_i$ contains all transitive tournaments on $\{1,\dots,n_i\}$.
        
        \item \label{label:tractCases2} There is $n\geq1$ such that the transitive tournament on $n$ vertices is the largest finite transitive tournament that is $\cF$-free and every tournament on strictly more vertices is not $\cF$-free. Moreover the standard representation of every $R_i$ contains all transitive tournaments on $\{1,\dots,n_i\}$. 
        
        \item \label{label:tractCases3} For every $n \leq m_\cF$, the set of $\cF$-free tournaments on $\{1,\dots,n\}$ as well as the (standard representations of the) relations $R_1,\dots R_m$ are preserved by the minority operation. 
        
        \item \label{label:tractCases4} For every $n \leq m_\cF$, the set of $\cF$-free tournaments on $\{1,\dots,n\}$ as well as the (standard representations of the) relations $R_1,\dots,R_m$ are preserved by the majority operation. 
    \end{enumerate}
Given the standard representation of every $R_i$ as well as the finite set of forbidden tournaments $\cF$, it can be checked algorithmically whether one of these cases applies.
\end{corollary}
\begin{proof}
    Denote the signature of $\bA$ by $\tau$. We first check that in every case of the list $\CSP(\bA)$ is in $P$.
    Tractability of cases~\ref{label:tractCases1}.~and~\ref{label:tractCases2}. follows by simply describing the class of finite $\tau$-structures that homomorphically map to $\bA$. 
    
    Namely, suppose we are in the situation of item~\ref{label:tractCases1}.~and let $\bX=(X,E^\bX,R_1^\bX,\dots,R_n^\bX)$ be a finite $\tau$-structure. Then $\bX$ homomorphically maps to $\bA$ if and only if $E^\bX$ contains no loops, i.e., a tuple $(x,x)$ for some $x \in X$, and all $R_i^\bX$ consist only of injective tuples.
    The latter condition is clearly necessary so it remains to show that it is also sufficient. Let $\widetilde{E^\bX}$ be the symmetric closure of the binary relation $E^\bX$and pick an acyclic orientation $(X,\to^X)$ of the graph $(X,E^{\tilde \bX})$, by restricting an arbitrary linear order of $X$ to the edge relation $\widetilde{E^\bX}$. Since every transitive tournament is $\cF$-free there is an embedding $f\colon (X,\to^X) \to \DF$, which also is a homomorphism $(X,\widetilde{E^\bX}) \to \HF$, and using the description of the $R_i$ it is not hard to see that $f$ even is a homomorphism $\bX \to \bA$.
    
    If we are in the situation of item~\ref{label:tractCases2}.~a finite $\tau$-structure $\bX$ homomorphically maps to $\bA$ if and only if it meets the same conditions as above in item~\ref{label:tractCases1}., and additionally the graph $(X,\widetilde{E^\bX})$, where $\widetilde{E^\bX}$ is the symmetric closure of $E^\bX$, does not contain an $n+1$-clique. The same reasoning as above works in this case. 

    Lastly, if the $\cF$-free tournaments on $\{1,\dots,n\}$ as well as all the $R_i$ are preserved by the minority (majority) operation, then Lemma~\ref{lem:buildingCanonicalFunction} stipulates the existence of a ternary $\DF$-canonical $f\in \Pol(\bA)$ that acts as a minority (majority) on the $\DF$-orbits $\{\xto,\to\}$. Letting $g$ be a ternary injection provided by Lemma~\ref{lem:theNdominatedfunctionG} we obtain
    \begin{align*}
        \tilde f(x,y,z) := f( g(x,y,z), g(y,z,x), g(z,x,y) ).
    \end{align*}
    The function $\tilde f$ is $\DF$-canonical an it is easily seen to act in cyclic fashion on the $\DF$-orbits of pairs $\{=,N,\xto,\to\}$. Hence the function $\tilde f$ is pseudo-cyclic modulo $\overline{\Aut(\DF)}$, see~\cite[Lemma 3]{BP-canonical}, and the reduction to the orbit-CSP~\cite{Bodirsky-Mottet} shows tractability of $\CSP(\bA)$. 

    Now we assume that $\CSP(\bA)$ is not NP-complete and show that one of the four cases holds. We first consider the situation that $\bA$ is not novel. If cases \ref{label:nonNovel1}.-\ref{label:nonNovel3}.~of Lemma~\ref{lem:nonNovelExplicit} apply we clearly are in one of the above cases. So assume that we are in the last case of said lemma, that is, $\bA'=(K_\omega,R_1'\dots,R_m')$ is a first-order reduct of the universal homogeneous tournament $\bT$ and it has a tractable CSP. By~\cite[Theorem 14]{MottetPinskerSmooth} there is a $\bT$-canonical ternary $f\in \Pol(\bA')$ that is pseudo-cyclic modulo $\overline{\Aut(\bT)}$. The function $f$ acts in cyclic fashion on the set of $\bT$-orbits of pairs $\{\xto,\to\}$ and by possibly replacing $f$ with $f\circ (f,f,f)$ it does so idempotently, i.e., either as a minority or as a majority operation. The same argument as in the proof of Lemma~\ref{lem:buildingCanonicalFunction} shows that all relations $R_1',\dots,R_m'$ are preserved by either the minority or majority operation. This immediately implies that all of $R_1,\dots,R_m$ are preserved by the minority or majority operation. 

    Lastly assume that $\bA$ is a novel structure. In this case Proposition~\ref{prop:MainThmForNovelStructures} guarantees the existence of a ternary $\DF$-canonical structure $f\in \Pol(\bA)$ that is pseudo-cyclic modulo $\overline{\Aut(\DF)}$. In particular $f$ acts cyclically on the set of $\DF$ orbits of edges in $\HF$, $\{\xto,\to\}$. Again, we may even assume that this action is idempotent, so $f$ either acts as a minority or as a majority and the desired result follows by invoking Lemma~\ref{lem:buildingCanonicalFunction}. 

    The criterion from the last two items can clearly be checked algorithmically. The criterion from the first two items can be checked algorithmically as the transitive tournament on $m_\cF$ vertices is $\cF$-free if and only if every transitive tournament is $\cF$-free. 
\end{proof}

\subsection{Examples}\label{sect:examples}
  
\begin{figure}[ht!]
\centering
\begin{tikzpicture}
  \begin{scope}[xshift = -4.5cm, scale = 1]
    \node [vertex] (1) at (-0.75,-0.75) {};
    \node [vertex] (2) at (0.75,-0.75) {};
    \node [vertex] (3) at (-0.75,0.75) {};
    \node [vertex] (4) at (0.75,0.75) {};
    \node (L1) at (0,-1.25) {$T_4$};
      
    \foreach \from/\to in {1/2, 1/3, 1/4, 2/3, 2/4, 3/4}     
    \draw [arc] (\from) to (\to);
  \end{scope}

  \begin{scope}[xshift = -1.5cm, scale = 1]
    \node [vertex] (1) at (-0.75,-0.75) {};
    \node [vertex] (2) at (0.75,-0.75) {};
    \node [vertex] (3) at (-0.75,0.75) {};
    \node [vertex] (4) at (0.75,0.75) {};
    \node (L1) at (0,-1.25) {$TC_4$};
      
    \foreach \from/\to in {1/2, 3/1, 1/4, 2/3, 2/4, 4/3}     
    \draw [arc] (\from) to (\to);
  \end{scope}

  \begin{scope}[xshift = 1.5cm, scale = 1]
    \node [vertex] (1) at (-0.75,-0.75) {};
    \node [vertex] (2) at (0.75,-0.75) {};
    \node [vertex] (3) at (-0.75,0.75) {};
    \node [vertex] (4) at (0.75,0.75) {};
    \node (L1) at (0,-1.25) {$C_3^{-}$};
      
    \foreach \from/\to in {1/2, 1/3, 1/4, 2/3, 4/2, 3/4}     
    \draw [arc] (\from) to (\to);
  \end{scope}

  \begin{scope}[xshift = 4.5cm, scale = 1]
    \node [vertex] (1) at (-0.75,-0.75) {};
    \node [vertex] (2) at (0.75,-0.75) {};
    \node [vertex] (3) at (-0.75,0.75) {};
    \node [vertex] (4) at (0.75,0.75) {};
    \node (L1) at (0,-1.25) {$C_3^{+}$};
      
    \foreach \from/\to in {2/1/, 3/1, 4/1, 2/3, 4/2, 3/4}     
    \draw [arc] (\from) to (\to);
  \end{scope}
\end{tikzpicture}
\caption{The tournaments on four vertices}
\label{fig:fourtournaments}
\end{figure}

We first exhibit a set of forbidden tournaments for which the orientation as well as the orientation-completion problem are tractable, while a generalised orientation problem is hard. 

Let $\cF$ consist of $T_4$ and $TC_4$ depicted in Figure~\ref{fig:fourtournaments}. 
In this case, the set $\cD$ consisting of $\cF$-free tournaments on $\{1,\dots,4\}$, i.e., those tournaments isomorphic to either $C_3^-$ or $C_3^+$, is preserved by the minority operation, as can be shown by an explicit computation. It even holds that $\mino$ acts on the isomorphism types of the $\cF$-free tournaments in this case, that is if exactly one or three of $D_1,D_2,D_3$ are isomorphic to $C_3^-$, then $\mino(D_1,D_2,D_3)$ is isomorphic to $C_3^-$ and the analogous statement holds for $C_3^+$.
Setting $\bA_1 = \HF$, and $\bA_2=(\HF,\to)$, templates that model the $\cF$-free orientation and the $\cF$-free orientation completion problem respectively, the third item of Corollary~\ref{cor:nonTechnicalVersion} implies that both computational problems are tractable. Now consider the relation $R\subseteq V^3$ consisting of all triples that induce the transitive tournament $T_3$ in $\DF$. 
Then $\CSP(\HF,R)$ models the computational problem where an input is a graph with some of its $3$-cliques labelled and the task consists of orienting the edges of this graph in a way that all of the labelled $3$-cliques are oriented transitively. 
Figure~\ref{fig:notPreservedByMinMaj} demonstrates that $R$ is neither preserved by the minority operation nor by the majority operation so Corollary~\ref{cor:nonTechnicalVersion} implies that the problem is NP-complete.

\begin{figure}
\centering
\begin{tikzpicture}

  \begin{scope}[xshift = -6cm]
    \node [vertex] (1) at (0,0) {};
    \node [vertex] (2) at (1,0) {};
    \node [vertex] (3) at (0.5,0.86602540378) {};
    \node (L1) at (0.5,-0.5) {$a \in R$};
      
    \foreach \from/\to in {1/2, 1/3, 2/3}     
    \draw [arc] (\from) to (\to);
  \end{scope}

  \begin{scope}[xshift = -3.5cm]
    \node [vertex] (1) at (0,0) {};
    \node [vertex] (2) at (1,0) {};
    \node [vertex] (3) at (0.5,0.86602540378) {};
    \node (L1) at (0.5,-0.5) {$b \in R$};
      
    \foreach \from/\to in {1/2, 3/1, 3/2}     
    \draw [arc] (\from) to (\to);
  \end{scope}

  \begin{scope}[xshift = -1cm]
    \node [vertex] (1) at (0,0) {};
    \node [vertex] (2) at (1,0) {};
    \node [vertex] (3) at (0.5,0.86602540378) {};
    \node (L1) at (0.5,-0.5) {$c \in R$};
      
    \foreach \from/\to in {2/1, 3/1, 2/3}     
    \draw [arc] (\from) to (\to);
  \end{scope}

  \begin{scope}[xshift = 2.5cm]
    \node [vertex] (1) at (0,0) {};
    \node [vertex] (2) at (1,0) {};
    \node [vertex] (3) at (0.5,0.86602540378) {};
    \node (L1) at (0.5,-0.5) {$\mino(a,b,c) \not\in R$};
      
    \foreach \from/\to in {2/1, 1/3, 3/2}     
    \draw [arc] (\from) to (\to);
  \end{scope}

  \begin{scope}[xshift = 6cm]
    \node [vertex] (1) at (0,0) {};
    \node [vertex] (2) at (1,0) {};
    \node [vertex] (3) at (0.5,0.86602540378) {};
    \node (L1) at (0.5,-0.5) {$\majo(a,b,c) \not \in R$};
      
    \foreach \from/\to in {1/2, 3/1, 2/3}     
    \draw [arc] (\from) to (\to);
  \end{scope}
  
\end{tikzpicture}
\caption{The relation $R$ is not closed under $\mino$ nor $\majo$}
\label{fig:notPreservedByMinMaj}
\end{figure}

As a second and last example we describe a generalised orientation problem where tractability is witnessed by a ternary $\DF$-canonical polymorphism that acts as a majority on $\{\xto,\to\}$, but no ternary canonical minority operation exists. This contrasts the result from~\cite{forbiddenTournaments}, which implies that tractability of the orientation and orientation completion problem, excluding the trivial cases from item~\ref{label:tractCases1}.~and~\ref{label:tractCases2}.~of Corollary~\ref{cor:nonTechnicalVersion}, is always accompanied by the existence of a ternary $\DF$-canonical canonical function which acts as a minority on $\{\xto,\to\}$. 
Simply take $\cF$ to be empty and let $R$ be the relation 
from Section~\ref{sect:minorityMajority}, which served as an example of a relation that was preserved by the majority operation but not by the minority operation. Clearly $\CSP(\HF,R)$ is not among the trivial cases of Corollary~\ref{cor:nonTechnicalVersion} and its tractability is witness by a canonical majority function, while no ternary canonical function acting as a minority on $\{\xto,\to\}$ can exist, by Lemma~\ref{lem:buildingCanonicalFunction}.  

\section{Outlook}\label{sect:outlook}

We hope that the approach  exhibited in the present paper provides a starting point for  complexity classifications of generalisations of the here discussed graph orientation problems. Such generalisations include the following.

\begin{itemize}
\item Graph orientation problems with forbidden tournaments where we are given a graph enhanced with additional local constraints which prescribe some labelled tuples of vertices of the given graph to belong to a fixed set of allowed directed graphs after orientation. Such problems can be modeled  as the CSP of structures of the form  $(\HF,R_1,\dots,R_m)$ where the $\DF$-definable relations $R_1,\dots,R_m$ are no longer required to induce tournaments in the graph $\DF$. This would, in particular, simultaneously generalize our results and the classification of Graph-SAT problems from~\cite{BodPin-Schaefer-both}.\label{1}
\item The $\cF$-free orientation problem for sets $\cF$ of oriented graphs which are  not necessarily tournaments. Using a theorem of Cherlin, Shelah, and Shi~\cite{CherlinShelahShi}, and its refinement by Hubi\v{c}ka and Ne\v{s}et\v{r}il~\cite{CSS-Ramsey}, this problem can be modelled as the  constraint satisfaction problem of a template enjoying similar properties to that modelling the situation where $\cF$ consists of tournaments, provided that $\cF$ is closed under homomorphic images. More precisely, the template would then be a first-order reduct of a finitely bounded homogeneous Ramsey structure which is an expansion of $(\DF,<)$. 
\item More general  edge-colouring problems as  captured by the logic GMSNP~\cite{bienvenu2014}, whose complexity classification belongs to the  most prominent open problems in infinite-domain constraint satisfaction.
\end{itemize}

While we hope that our methods generalise smoothly to computational problems in the first item, the second and third item appear to pose greater challenges. One obstacle is the lack of a known description of the automorphism groups of the first-order structures parametrising the respective computational problems, i.e., of an analogue of the result in~\cite{ClassificationOfReductsOfDF}. Although the full strength of that result  does not seem to be necessary, a statement similar to Corollary~\ref{NormalFormAutomorphisms} providing  ``sufficient  control'' over the action of automorphisms would still be required.
Furthermore, in contrast to the present paper, where we could in many cases restrict ourselves to studying various clones acting on the two-element set $\{\xto,\to\}$ of possible orientations of an edge, more general cases will involve clones acting on larger finite sets. For example, in edge-colouring problems, we may need to consider the actions on the possible colourings of an edge by $n$ many colours $\{c_1,\dots,c_n\}$. On the one hand this leads to more involved combinatorics in some cases, e.g., in Lemmas~\ref{gettingMajority} and~\ref{lem:gettingMajFromCommutativeFunction}, while on the other hand it precludes the direct application of Post's classification of clones on a two-element domain~\cite{Post}, exploited in Lemma~\ref{DFcanimpliesHFcan}.
 
Concretely, the problem of colouring a graph with $n$ colours avoiding a list of $n$-coloured cliques could serve as a natural first step in advancing the results from the present paper. The template modelling such a problem is still built from a free amalgamation class, just as in this paper, making an analogous result to~\cite{ClassificationOfReductsOfDF}, or rather, a sufficient corollary seem  within reach. This would allow one to focus on the core challenge of dealing with more than two “relevant” orbits.

\bibliographystyle{plain}
\bibliography{references}

\end{document}